\documentclass[11pt]{article}%
\usepackage{amsfonts}
\usepackage{graphicx}
\usepackage{amsmath}
\usepackage{amssymb}
\usepackage{dsfont}
\usepackage[a4paper,bottom=3.5cm,top=2cm,left=2.5cm,right=2.5cm]{geometry}%
\providecommand{\U}[1]{\protect\rule{.1in}{.1in}}
\newcommand\beq{\begin{equation}}
\newcommand\eeq{\end{equation}}

\newtheorem{theorem}{Theorem}

\newtheorem{corollary}[theorem]{Corollary}

\newtheorem{lemma}[theorem]{Lemma}

\newtheorem{proposition}[theorem]{Proposition}
\newtheorem{remark}[theorem]{Remark}

\begin{document}

\begin{center}
{\Large \textbf{On the limiting spectral distribution for a large class of
symmetric random matrices with correlated entries.}} \medskip

Marwa Banna$^{a}$, Florence Merlev\`{e}de$^{b}$, Magda Peligrad$^{b}%
$\footnote{Supported in part by a Charles Phelps Taft Memorial Fund grant, and
the NSF grant DMS-1208237.}
\end{center}

$^{a \, b}$ Universit\'{e} Paris Est, LAMA (UMR 8050), UPEMLV, CNRS, UPEC, 5
Boulevard Descartes, 77 454 Marne La Vall\'{e}e, France.

E-mail: marwa.banna@u-pem.fr; florence.merlevede@u-pem.fr

$^{c}$ Department of Mathematical Sciences, University of Cincinnati, PO Box
210025, Cincinnati, Oh 45221-0025, USA.

Email: peligrm@ucmail.uc.edu

\bigskip



\begin{center}
\bigskip

\bigskip\textbf{Abstract}
\end{center}

For symmetric random matrices with correlated entries, which are functions of
independent random variables, we show that the asymptotic behavior of the
empirical eigenvalue distribution can be obtained by analyzing a Gaussian
matrix with the same covariance structure. This class contains both cases of
short and long range dependent random fields. The technique is based on a
blend of blocking procedure and Lindeberg's method. This method leads to a
variety of interesting asymptotic results for matrices with dependent entries,
including applications to linear processes as well as nonlinear Volterra-type
processes entries.

\section{Introduction}

The limiting spectral distribution for symmetric matrices with correlated
entries received a lot of attention in the last two decades. The starting
point is deep results for symmetric matrices with correlated Gaussian entries
by Khorunzhy and Pastur \cite{KP}, Boutet de Monvel \textit{et al} \cite{BKV},
Boutet de Monvel and Khorunzhy \cite{BK}, Chakrabarty \textit{et al}
\cite{chak} among others. On the other hand there is a sustained effort for
studying linear filters of independent random variables as entries of a
matrix. For instance, Anderson and Zeitouni \cite{AZ} considered symmetric
matrices with entries that are linear processes of finite range of independent
random variables. Hachem \textit{et al} \cite{HLN} considered large sample
covariance matrices whose entries are modeled by a short memory linear process
of infinite range with independent Gaussian innovations. Bai and Zhou
\cite{BZ}, Yao \cite{Yao}, Banna and Merlev\`{e}de \cite{BM} and Merlev\`{e}de
and Peligrad \cite{MP}, among others, treated large covariance matrices based
on an independent structure between columns and correlated random variables in rows.

In this paper we consider symmetric random matrices whose entries are
functions of independent and identically distributed (i.i.d.) real-valued
random variables. Such kind of processes provide a very general framework for
stationary ergodic random fields. Our main goal is to reduce the study of the
limiting spectral distribution to the same problem for a Gaussian matrix
having the same covariance structure as the underlying process. In this way we
prove the universality and we are able to formulate various limiting results
for large classes of matrices. We also treat large sample covariance matrices
with correlated entries, known under the name of Gram matrices. Our proofs are
based on the large-small block arguments, a method which, in one dimensional
setting, is going back to Bernstein. Then, we apply a variant of the so-called
Lindeberg method, namely we develop a block Lindeberg method, where we replace
at one time a big block of random variables with a Gaussian one with the same
covariance structure. Lindeberg method is popular with these type of problems.
Replacing only one variable at one time with a Gaussian one, Chatterjee
\cite{Ct} treated random matrices with exchangeable entries.

Our paper is organized in the following way. Section \ref{SectMR} contains the
main results for symmetric random matrices and sample covariance matrices. As
an intermediate step we also treat matrices based on K-dependent random
fields, results that have interest in themselves (see Theorem
\ref{thmmdependent} in Section \ref{sectionKdependent}). In Section
\ref{Sectexamples}, we give applications to matrices with entries which are
either linear random fields or nonlinear random fields as Volterra-type
processes. The main proofs are included in Section \ref{Sectproofs}. In
Section \ref{Sectappendix} we prove a concentration of spectral measure
inequality for a row-wise K-dependent random matrix and we also mention some
of the technical results used in the paper.

Here are some notations used all along the paper. The notation $[x]$ is used
to denote the integer part of any real $x$. For any positive integers $a,b$,
the notation $\mathbf{0}_{a,b}$ means a matrix with $0$ entries of size
$a\times b$, whereas the notation $\mathbf{0}_{a}$ means a row vector of size
$a$. For a matrix $A$, we denote by $A^{T}$ its transpose matrix and by
$\mathrm{Tr}(A)$ its trace. We shall use the notation $\Vert X\Vert_{r}$ for
the ${\mathbb{L}}^{r}$-norm ($r\geq1$) of a real valued random variable $X$.

For any square matrix $A$ of order $n$ with only real eigenvalues $\lambda
_{1}\leq\dots\leq\lambda_{n}$, its spectral empirical measure and its spectral
distribution function are respectively defined by
\[
\nu_{A}=\frac{1}{n}\sum_{i=1}^{n}\delta_{\lambda_{i}}\ \text{ and }\ F_{n}%
^{A}(x)=\frac{1}{n}\sum_{k=1}^{n}\mathbf{1}_{\{\lambda_{k}\leq x\}}\,.
\]
The Stieltjes transform of $F^{A}$ is given by
\[
S_{A}(z)=\int\frac{1}{x-z}dF^{A}(x)=\frac{1}{n}\mathrm{Tr}(A-z{\mathbf{I}}%
_{n})^{-1}\,,
\]
where $z=u+\mathrm{i}v\in\mathbb{C}^{+}$ (the set of complex numbers with
positive imaginary part), and $\mathbf{{I}}_{n}$ is the identity matrix of
order $n$.

The L\'evy distance between two distribution functions $F$ and $G$ is defined
by
\[
L(F,G)=\inf\{\varepsilon>0\ :\ F(x-\varepsilon)-\varepsilon\leq G(x)\leq
F(x+\varepsilon)+\varepsilon\}\,.
\]
It is well-known that a sequence of distribution functions $F_{n}(x)$
converges to a distribution function $F(x)$ at all continuity points $x$ of
$F$ if and only if $L(F_{n},F)\rightarrow0$.

\medskip

\section{Main results}

\label{SectMR}

\subsection{On the limiting distribution for a large class of symmetric
matrices with correlated entries}

Let $(X_{k,\ell})_{(k,\ell)\in{\mathbb{Z}}^{2}}$ be an array of real-valued
random variables, and consider its associated symmetric random matrix
$\mathbf{X}_{n}$ of order $n$ defined by
\begin{equation}
\text{$\big (\mathbf{X}_{n}\big )_{i,j}=X_{i,j}$ if $1\leq j\leq i\leq n$ and
$\big (\mathbf{X}_{n}\big )_{i,j}=X_{j,i}$ if $1\leq i<j\leq n$}\,.
\label{defbfX}%
\end{equation}
Define then
\begin{equation}
{\mathbb{X}}_{n}:=n^{-1/2}\mathbf{X}_{n}\,. \label{defbbX}%
\end{equation}
The aim of this section is to study the limiting spectral empirical
distribution function of the symmetric matrix ${\mathbb{X}}_{n}$ defined by
\eqref{defbbX} when the process $(X_{k,\ell})_{(k,\ell)\in{\mathbb{Z}}^{2}}$
has the following dependence structure: for any $(k,\ell)\in{\mathbb{Z}}^{2}%
$,
\begin{equation}
X_{k,\ell}=g(\xi_{k-i,\ell-j}\ ;\ (i,j)\in{\mathbb{Z}}^{2})\,, \label{defYkl}%
\end{equation}
where $(\xi_{i,j})_{(i,j)\in{\mathbb{Z}}^{2}}$ is an array of i.i.d.
real-valued random variables given on a common probability space
$(\Omega,\mathcal{K},\mathbb{P})$, and $g$ is a measurable function from
${\mathbb{R}}^{{\mathbb{Z}}^{2}}$ to ${\mathbb{R}}$ such that ${\mathbb{E}%
}(X_{0,0})=0$ and $\Vert X_{0,0}\Vert_{2}<\infty$. A representation as
\eqref{defYkl} includes linear as well as many widely used nonlinear random
fields models as special cases.

\smallskip Our Theorem \ref{thmUni} below shows a universality scheme for the
random matrix ${\mathbb{X}}_{n}$ as soon as the entries of the symmetric
matrix $\sqrt n {\mathbb{X}}_{n}$ have the dependence structure
\eqref{defYkl}. It is noteworthy to indicate that this result does not require
rate of convergence to zero of the correlation between the entries.

\begin{theorem}
\label{thmUni} Let $(X_{k,\ell})_{(k,\ell)\in{\mathbb{Z}}^{2}}$ be a
real-valued stationary random field given by \eqref{defYkl}. Define the
symmetric matrix ${\mathbb{X}}_{n}$ by \eqref{defbbX}. Let $(G_{k,\ell
})_{(k,\ell)\in{\mathbb{Z}}^{2}}$ be a real-valued centered Gaussian random
field, with covariance function given by
\begin{equation}
{\mathbb{E}}(G_{k,\ell}G_{i,j})={\mathbb{E}}(X_{k,\ell}X_{i,j})\text{ for any
$(k,\ell)$ and $(i,j)$ in ${\mathbb{Z}}^{2}$}\,. \label{egacovfunction}%
\end{equation}
Let $\mathbf{G}_{n}$ be the symmetric random matrix defined by
$\big (\mathbf{G}_{n}\big )_{i,j}=G_{i,j}$ if $1\leq j\leq i\leq n$ and
$\big (\mathbf{G}_{n}\big )_{i,j}=G_{j,i}$ if $1\leq i<j\leq n$. Denote
${\mathbb{G}}_{n}=\frac{1}{\sqrt{n}}{\mathbf{G}}_{n}$. Then, for any
$z\in\mathbb{C}^{+}$,
\[
\lim_{n\rightarrow\infty}\big |S_{{\mathbb{X}}_{n}}(z)-{\mathbb{E}%
}\big (S_{{\mathbb{G}}_{n}}(z)\big )\big |=0\ \text{ almost surely.}%
\]

\end{theorem}

Theorem \ref{thmUni} is important since it shows that the study of the
limiting spectral distribution function of a symmetric matrix whose entries
are functions of i.i.d. random variables can be reduced to studying the same
problem as for a Gaussian matrix with the same covariance structure. The
following corollary is a direct consequence of our Theorem \ref{thmUni}
together with Theorem B.9 in Bai-Silverstein \cite{BS} (see also the arguments
on page 38 in \cite{BS}, based on Vitali's convergence theorem).

\begin{corollary}
\label{coruni} Assume that ${\mathbb{X}}_{n}$ and $\mathbb{G}_{n}$ are as in
Theorem \ref{thmUni}. Furthermore, assume there exists a distribution function
$F$ such that
\[
{\mathbb{E}}\big (F^{{\mathbb{G}}_{n}}(t)\big )\rightarrow F(t)\text{ for all
continuity points }t\in{\mathbb{R}}\text{ of }F\,.
\]
Then
\begin{equation}
{\mathbb{P}}(L(F^{{\mathbb{X}}_{n}(\omega)},F)\rightarrow0)=1\, .
\label{conv in distribution}%
\end{equation}

\end{corollary}

For instance, Corollary \ref{coruni} above combined with the proof of Theorem
2 in Khorunzhy and Pastur \cite{KP} concerning the asymptotic spectral
behavior of certain ensembles with correlated Gaussian entries (see also
Theorem 17.2.1 in \cite{PS}), gives the following:

\begin{theorem}
\label{SCL1} Let $(X_{k,\ell})_{(k,\ell)\in{\mathbb{Z}}^{2}}$ be a real-valued
stationary random field given by \eqref{defYkl}. Define the symmetric matrix
${\mathbb{X}}_{n}$ by \eqref{defbbX}. For any $(k,\ell)\in{\mathbb{Z}}^{2}$,
let $\gamma_{k,\ell}={\mathbb{E}}(X_{0,0}X_{k,\ell})$. Assume that
\begin{equation}
\sum_{k,\ell\in\mathbb{Z}} |\gamma_{k,\ell} | < \infty\,
,\label{condabssumcov}%
\end{equation}
and that the following holds: for any $(k,\ell)\in{\mathbb{Z}}^{2}$,
\begin{equation}
\gamma_{k,\ell}=\gamma_{\ell,k}\, , \label{condgammaY}%
\end{equation}
Then (\ref{conv in distribution}) holds, where $F$ is a nonrandom distribution
function whose Stieltjes transform $S(z)$ is uniquely defined by the
relations:
\begin{equation}
S(z)=\int_{0}^{1}h(x,z)dx\,, \label{defSlim}%
\end{equation}
where $h(x,z)$ is a solution to the equation
\begin{equation}
h(x,z)=\Big (-z+\int_{0}^{1}f(x,y)h(y,z)dy\Big )^{-1} \ \text{ with }
\ f(x,y)= \sum_{k,j \in\mathbb{Z}}\gamma_{k,j} {\mathrm{e}}^{-2\pi
\mathrm{i}(kx+jy)} \, . \label{defSlim2}%
\end{equation}

\end{theorem}

Equation \eqref{defSlim2} is uniquely solvable in the class ${\mathcal{F}}$ of
functions $h(x,z)$ with domain $(x,z)\in\lbrack0,1]\otimes{\mathbb{C}%
}\backslash{\mathbb{R}}$, which are analytic with respect to $z$ for each
fixed $x$, continuous with respect to $x$ for each fixed $z$ and satisfying
the conditions: $\lim_{v\rightarrow\infty}v\operatorname{Im}h(x,\mathrm{i}%
v)\leq1$ and $\operatorname{Im}(z)\operatorname{Im}h(x,z)>0$.

\begin{remark}
\label{remthG}If condition \eqref{condgammaY} of Theorem \ref{SCL1} is
replaced by: $\gamma_{\ell,k}=V(\ell)V(k)$ where $V$ is an even function, then
its conclusion can be given in the following alternative way: the convergence
(\ref{conv in distribution}) holds where $F$ is a nonrandom distribution
function whose Stieltjes transform $S(z)$ is given by the relation
\[
S(z)=\int_{0}^{\infty}\frac{d\upsilon(\lambda)}{-z-\lambda h(z)}%
\]
where $\upsilon(t)=\lambda\{x\in\lbrack0,1];f(x)<t\}$, $\lambda$ is the
Lebesgue measure, for $x\in\lbrack0,1]$, $f(x)=\sum_{k\in{\mathbb{Z}}%
}V(k)\mathrm{e}^{2\pi\mathrm{i}kx}$ and $h(z)$ is solution to the equation
\[
h(z)=\int_{0}^{\infty}\frac{\lambda d\upsilon(\lambda)}{-z-\lambda
h(z)}\text{,\ \ }z\in{\mathbb{C}}\backslash{\mathbb{R}}\,.
\]
This equation is uniquely solvable in the class of analytic functions in
${\mathbb{C}}\backslash{\mathbb{R}}$ satisfying the conditions: $\lim
_{x\rightarrow\infty}xh(\mathrm{i}x)<\infty$ and $\operatorname{Im}%
\big (h(z)\big )\operatorname{Im}(z)>0$ for $z\in{\mathbb{C}}\backslash
{\mathbb{R}}$. (See Boutet de Monvel and Khorunzhy \cite{BK}).
\end{remark}

\subsection{On the limiting distribution for large Gram (sample covariance)
matrices with correlated entries}

Adapting the proof of Theorem \ref{thmUni}, we can also obtain a universality
scheme for large sample covariance matrices associated with a process
$(X_{k,\ell})_{(k,\ell)\in{\mathbb{Z}}^{2}}$ having the representation
\eqref{defYkl}. So, all along this section $(X_{k,\ell})_{(k,\ell
)\in{\mathbb{Z}}^{2}}$ is assumed to be a random field having the
representation \eqref{defYkl}. To define the Gram matrices associated with
this random field, we consider two positive integers $N$ and $p$, and the
$N\times p$ matrix
\begin{equation}
{\mathcal{X}}_{N,p}=\big (X_{i,j}\big )_{1\leq i\leq N,1\leq j\leq p}\, .
\label{defcalX}%
\end{equation}
Define now the symmetric matrix ${\mathbb{B}}_{N}$ of order $N$ by
\begin{equation}
{\mathbb{B}}_{N}=\frac{1}{p}{\mathcal{X}}_{N,p}{\mathcal{X}}_{N,p}^{T}: =
\frac{1}{p} \sum_{k=1}^{p} {\mathbf{r}}_{k} {\mathbf{r}}_{k}^{T} \,
,\label{defBN}%
\end{equation}
where ${\mathbf{r}}_{k} = ( X_{1,k}, \cdots, X_{N,k} )^{T}$ is the $k$-th
column of ${\mathcal{X}}_{N,p}$.

The matrix ${\mathbb{B}}_{N}$ is usually referred to as the sample covariance
matrix associated with the process $(X_{k,\ell})_{(k,\ell)\in{\mathbb{Z}}^{2}%
}$. It is also known under the name of Gram random matrix.

\begin{theorem}
\label{thmainGram}Let ${\mathbb{B}}_{N}$ be defined by (\ref{defBN}) and let
${\mathbb{H}}_{N}= \frac{1}{p}{\mathcal{G}}_{N,p}{\mathcal{G}}_{N,p}^{T} $ be
the Gram matrix associated with a real-valued centered Gaussian random field
$(G_{k,\ell})_{(k,\ell)\in{\mathbb{Z}}^{2}}$, with covariance function given
by (\ref{egacovfunction}). Then, provided that $N,p\rightarrow\infty$ such
that $N/p\rightarrow c \in(0,\infty)$, for any $z\in\mathbb{C}^{+}$,
\begin{equation}
\lim_{n\rightarrow\infty}\big |S_{{\mathbb{B}}_{N}}(z)-{\mathbb{E}%
}(S_{{\mathbb{H}}_{N}}(z)\big )\big |=0\ \text{ almost surely.}
\label{thmainGramcomp}%
\end{equation}
Therefore, if $N,p\rightarrow\infty$ such that $N/p\rightarrow c \in
(0,\infty)$ and if there exists a distribution function $F$ such that
\[
{\mathbb{E}}\big (F^{{\mathbb{H}}_{N}}(t)\big )\rightarrow F(t)\text{ for all
continuity points }t\in{\mathbb{R}}\text{ of }F \,
\]
then
\begin{equation}
\label{conv in distribution gram}{\mathbb{P}}(L(F^{{\mathbb{B}}_{N}(\omega
)},F)\rightarrow0)=1\, .
\end{equation}

\end{theorem}

Theorem \ref{thmainGram} together with Theorem 2.1 in Boutet de Monvel
\textit{et al} \cite{BKV} allow then to derive the limiting spectral
distribution of large sample covariance matrices associated with a process
$(X_{k,\ell})_{(k,\ell)\in{\mathbb{Z}}^{2}}$ having the representation
\eqref{defYkl} and satisfying a short range dependence condition.

\begin{theorem}
\label{MPthm} Let $(X_{k,\ell})_{(k,\ell)\in{\mathbb{Z}}^{2}}$ be a
real-valued stationary random field given by \eqref{defYkl}. Assume that
\eqref{condabssumcov} holds. Then, provided that $N,p\rightarrow\infty$ such
that $N/p\rightarrow c \in(0,\infty)$, ${\mathbb{P}}(L(F^{{\mathbb{B}}%
_{N}(\omega)},F)\rightarrow0)=1$ where $F$ is a nonrandom distribution
function whose Stieltjes transform $S(z)$, $z\in{\mathbb{C}}^{+}$ is uniquely
defined by the relations:
\[
S(z)=\int_{0}^{1}h(x,z)dx\,,
\]
where $h(x,z)$ is a solution to the equation
\begin{equation}
h(x,z)=\Big (-z+\int_{0}^{1}\frac{f(x,s)}{1+ c \int_{0}^{1}f(u,s)h(u,z)du}%
ds\Big )^{-1}\,, \label{eqg}%
\end{equation}
with $f(x,y)$ given in \eqref{defSlim2}.
\end{theorem}

Equation \eqref{eqg} is uniquely solvable in the class ${\mathcal{F}}$ of
functions $h(x,z)$ as described after the statement of Theorem \eqref{SCL1}.

\medskip

We refer to the paper by Boutet de Monvel \textit{et al} \cite{BKV} regarding
discussions on the smoothness and boundedness of the limiting density of
states. Note that condition \eqref{condabssumcov} is required in the statement
of Theorem \ref{MPthm} only because all the estimates in the proof of Theorem
2.1 in \cite{BKV} require this condition. However using arguments as developed
in the paper by Chakrabarty \textit{et al} \cite{chak}, it can be proved that
if the process $(X_{k,\ell})_{(k,\ell)\in{\mathbb{Z}}^{2}}$ admits a spectral
density then there exists a nonrandom distribution function $F$ such that
${\mathbb{P}}(L(F^{{\mathbb{B}}_{N}(\omega)},F)\rightarrow0)=1$ (if
$N/p\rightarrow c \in(0,\infty)$). Unfortunately the arguments developed in
\cite{chak} do not allow, in general, to exhibit the limiting equation
\eqref{eqg} which gives a lot of information on the limiting spectral
distribution. Notice however that if we add the assumption that the lines
(resp. the columns) of ${\mathcal{X}}_{N,p}$ are non correlated (corresponding to the semantically (resp. spatially) "patterns" studied in
Section 3 of \cite{BKV}), condition \eqref{condabssumcov} is not needed to
exhibit the limiting equation of the Stieltjes transform. Indeed, in this
situation, the lines (resp. the columns) of ${\mathcal{G}}_{N,p}$ become then
independent and the result of Merlev\`ede and Peligrad \cite{MP} about the
limiting spectral distribution of Gram random matrices associated to
independent copies of a stationary process applies. Proving, however, Theorem
\ref{MPthm} in its full generality and without requiring condition
\eqref{condabssumcov} to hold, remains an open question.

\section{Examples}

\label{Sectexamples} All along this section, $(\xi_{k,\ell})_{(k,\ell
)\in{\mathbb{Z}}^{2}}$ will designate a double indexed sequence of i.i.d.
real-valued random variables defined on a common probability space, centered
and in ${\mathbb{L}}^{2}$.

\subsection{Linear processes}

Let $(a_{k,\ell})_{(k,\ell)\in{\mathbb{Z}}^{2}}$ be a double indexed sequence
of numbers such that
\begin{equation}
\sum_{k,\ell\in{\mathbb{Z}}}|a_{k,\ell}|<\infty\,.\label{shortmem}%
\end{equation}
Let then $(X_{i,j})_{(i,j)\in{\mathbb{Z}}^{2}}$ be the linear random field in
${\mathbb{L}}^{2}$ defined by: for any $(i,j)\in
{\mathbb{Z}}^{2}$,
\begin{equation}
X_{i,j}=\sum_{k,\ell\in{\mathbb{Z}}}a_{k,\ell}\xi_{k+i,\ell+j}%
\,.\label{linproc}%
\end{equation}
Corollary \ref{coruni} (resp. Theorem \ref{thmainGram}) then applies to the
matrix $\mathbb{X}_{n}$ (resp. $\mathbb{B}_{N}$) associated with the linear
random field $(X_{i,j})_{(i,j)\in{\mathbb{Z}}^{2}}$ given in \eqref{defbbX}.

\medskip

In case of short dependence, based on our Theorem \ref{MPthm}, we can describe
the limit of the empirical spectral distribution of the Gram matrix associated
with a linear random field.

\begin{corollary}
Assume that $X_{i,j}$ is defined by (\ref{linproc}) and condition
(\ref{shortmem}) is satisfied. Let $N$ and $p$ be positive integers, such that
$N,p\rightarrow\infty$, $N/p\rightarrow c \in(0, \infty) $. Let ${\mathcal{X}%
}_{N,p}=\big (X_{i,j}\big )_{1\leq i\leq N,1\leq j\leq p}$ and ${\mathbb{B}%
}_{N}=N^{-1}{\mathcal{X}}_{N,p}{\mathcal{X}}_{N,p}^{T}$. Then the convergence
(\ref{conv in distribution gram}) holds for $F^{{\mathbb{B}}_{N}}$, where $F$
is a nonrandom distribution function whose Stieltjes transform satisfies the
relations given in Theorem \ref{MPthm} with $\gamma_{k,j}=\Vert\xi_{0,0}%
\Vert_{2}^{2}\sum_{u,v\in{\mathbb{Z}}}a_{u,v}a_{u+k,v+j}$.
\end{corollary}

Concerning now the Wigner-type matrix $\mathbb{X}_{n}$, by using Remark
\ref{remthG}, we obtain the following corollary, describing the limit in a
particular case.

\begin{corollary}
Let $(a_{n})_{n\in{\mathbb{Z}}}$ be a sequence of numbers such that
\[
\sum_{k\in{\mathbb{Z}}}|a_{k}|<\infty.
\]
Define $X_{i,j}=\sum_{k,\ell\in{\mathbb{Z}}}a_{k}a_{\ell}\xi_{k+i,\ell+j}$ for
any $(i,j)\in{\mathbb{Z}}^{2}$. Consider the symmetric matrix $\mathbb{X}_{n}$
associated with $(X_{i,j})_{(i,j)\in{\mathbb{Z}}^{2}}$ and defined by
\eqref{defbbX}. Then (\ref{conv in distribution}) holds, where $F$ is a
nonrandom distribution function whose Stieltjes transform satisfies the
relation given in Remark \ref{remthG} with $f(x)=\Vert\xi_{0,0}\Vert_{2}%
\sum_{k\in{\mathbb{Z}}}\sum_{j\in{\mathbb{Z}}}a_{j}a_{j+k}\mathrm{e}%
^{2\pi\mathrm{i}kx}$.
\end{corollary}

\bigskip

\subsection{Volterra-type processes}

Other classes of stationary random fields having the representation
\eqref{defYkl} are Volterra-type processes which play an important role in the
nonlinear system theory. For any ${\mathbf{k}}=(k_{1},k_{2})\in{\mathbb{Z}%
}^{2}$, define a second-order Volterra expansion as follows:
\begin{equation}
X_{{\mathbf{k}}}=\sum_{{\mathbf{u}}\in{\mathbb{Z}}^{2}}a_{{\mathbf{u}}}%
\xi_{{\mathbf{k}}-{\mathbf{u}}}+\sum_{{\mathbf{u}},{\mathbf{v}}\in{\mathbb{Z}%
}^{2}}b_{{\mathbf{u}},{\mathbf{v}}}\xi_{{\mathbf{k}}-{\mathbf{u}}}%
\xi_{{\mathbf{k}}-{\mathbf{v}}}\,, \label{Volterra}%
\end{equation}
where $a_{{\mathbf{u}}}$ and $b_{{\mathbf{u}},{\mathbf{v}}}$ are real numbers
satisfying
\begin{equation}
\label{condcoeffvolt}b_{{\mathbf{u}},{\mathbf{v}}}=0\text{ if }{\mathbf{u}%
}={\mathbf{v}}\,,\,\sum_{{\mathbf{u}}\in{\mathbb{Z}}^{2}}a_{{\mathbf{u}}}%
^{2}<\infty\ \text{ and }\ \sum_{{\mathbf{u}},{\mathbf{v}}\in{\mathbb{Z}}^{2}%
}b_{{\mathbf{u}},{\mathbf{v}}}^{2}<\infty\,.
\end{equation}
Under the above conditions, the random field $X_{{\mathbf{k}}}$ exists, is
centered and in ${\mathbb{L}}^{2}$, and Corollary \ref{coruni} (resp. Theorem
\ref{thmainGram}) applies to the matrix $\mathbb{X}_{n}$ (resp. $\mathbb{B}%
_{N}$) associated with the Volterra-type random field. Further generalization
to arbitrary finite order Volterra expansion is straightforward.

If we reinforced condition \eqref{condcoeffvolt}, we derive the following
result concerning the limit of the empirical spectral distribution of the Gram
matrix associated with the Volterra-type process:

\begin{corollary}
\label{propVolt1} Assume that $(X_{{\mathbf{k}}})_{{\mathbf{k}}\in{\mathbb{Z}%
}^{2}}$ is defined by (\ref{Volterra}) and that the following additional
condition is assumed:
\begin{equation}
\sum_{{\mathbf{u}}\in{\mathbb{Z}}^{2}}|a_{{\mathbf{u}}}|<\infty\ ,\ \sum
_{{\mathbf{v}}\in{\mathbb{Z}}^{2}}\Big (\sum_{{\mathbf{u}}\in{\mathbb{Z}}^{2}%
}b_{{\mathbf{u}},{\mathbf{v}}}^{2}\Big )^{1/2}<\infty\ \text{ and }%
\ \sum_{{\mathbf{v}}\in{\mathbb{Z}}^{2}}\Big (\sum_{{\mathbf{u}}\in
{\mathbb{Z}}^{2}}b_{{\mathbf{v}},{\mathbf{u}}}^{2}\Big )^{1/2}<\infty\,.
\label{condVolt1}%
\end{equation}
Let $N$ and $p$ be positive integers, such that $N,p\rightarrow\infty$,
$N/p\rightarrow c \in(0, \infty)$. Let ${\mathcal{X}}_{N,p}=\big (X_{i,j}%
\big )_{1\leq i\leq N,1\leq j\leq p}$ and ${\mathbb{B}}_{N}=N^{-1}%
{\mathcal{X}}_{N,p}{\mathcal{X}}_{N,p}^{T}$. Then
(\ref{conv in distribution gram}) holds for $F^{{\mathbb{B}}_{N}}$, where $F$
is a nonrandom distribution function whose Stieltjes transform satisfies the
relations given in Theorem \ref{MPthm} with
\begin{equation}
\gamma_{{\mathbf{k}}}=\Vert\xi_{0,0}\Vert_{2}^{2}\sum_{{\mathbf{u}}%
\in{\mathbb{Z}}^{2}}a_{{\mathbf{u}}}a_{{\mathbf{u}}+{\mathbf{k}}}+\Vert
\xi_{0,0}\Vert_{2}^{4}\sum_{{\mathbf{u}},{\mathbf{v}}\in{\mathbb{Z}}^{2}%
}b_{{\mathbf{u}},{\mathbf{v}}}(b_{{\mathbf{u}}+{\mathbf{k}},{\mathbf{v}%
}+{\mathbf{k}}}+b_{{\mathbf{v}}+{\mathbf{k}},{\mathbf{u}}+{\mathbf{k}}%
})\ \text{ for any ${\mathbf{k}}\in{\mathbb{Z}}^{2}$}\,. \label{defgammakvolt}%
\end{equation}

\end{corollary}

If we impose additional symmetric conditions to the coefficients $a_{
{\mathbf{u}}}$ and $b_{ {\mathbf{u}}, {\mathbf{v}}}$ defining the Volterra
random field \eqref{Volterra}, we can derive the limiting spectral
distribution of its associated symmetric matrix ${\mathbb{X}}_{n}$ defined by
\eqref{defbbX}. Indeed if for any ${\mathbf{u}}=(u_{1},u_{2})$ and
${\mathbf{v}}=(v_{1},v_{2})$ in ${\mathbb{Z}}^{2}$,
\begin{equation}
\label{condabVolt}a_{ {\mathbf{u}}}= a_{u_{1}}a_{u_{2}} \, , \, b_{
{\mathbf{u}}, {\mathbf{v}}}=b_{u_{1},v_{1}}b_{u_{2},v_{2}} \, ,
\end{equation}
where the $a_{i}$ and $b_{i,j}$ are real numbers satisfying
\begin{equation}
\label{condVolt2}\text{$b_{i,j}=0$ if $i=j$} \, , \, \sum_{i \in{\mathbb{Z}}}|
a_{ i}| < \infty\, \text{ and } \, \sum_{(i,j) \in{\mathbb{Z}}^{2}} |b_{i,j}|<
\infty\, ,
\end{equation}
then $(\gamma_{k,\ell})$ satisfies \eqref{condabssumcov} and
\eqref{condgammaY} . Hence, an application of Theorem \ref{SCL1} leads to the
following result.

\begin{corollary}
Assume that $(X_{{\mathbf{k}}})_{{\mathbf{k}}\in{\mathbb{Z}}^{2}}$ is defined
by (\ref{Volterra}) and that conditions \eqref{condabVolt} and
\eqref{condVolt2} are satisfied. Define the symmetric matrix ${\mathbb{X}}%
_{n}$ by \eqref{defbbX}. Then (\ref{conv in distribution}) holds, where $F$ is
a nonrandom distribution function whose Stieltjes transform is uniquely
defined by the relations given in Theorem \ref{SCL1} with $\gamma
_{s,t}=A(s)A(t)+B_{1}(s)B_{1}(t)+B_{2}(s)B_{2}(t)$ with $A(t)=\Vert\xi
_{0,0}\Vert_{2}\sum_{i\in{\mathbb{Z}}}a_{i}a_{i+t}$, $B_{1}(t)=\Vert\xi
_{0,0}\Vert_{2}^{2}\sum_{(i,r)\in{\mathbb{Z}}^{2}}b_{i,r}b_{i+t,r+t}$ and
$B_{2}(t)=\Vert\xi_{0,0}\Vert_{2}^{2}\sum_{(i,r)\in{\mathbb{Z}}^{2}}%
b_{i,r}b_{r+t,i+t}$.
\end{corollary}

\section{Proofs of the main results}

\label{Sectproofs}

The proof of Theorem \ref{thmUni} being based on an approximation of the
underlying symmetric matrix by a symmetric matrix with entries that are
$2m$-dependent (for $m$ a sequence of integers tending to infinity after $n$),
we shall first prove a universality scheme for symmetric matrices with
$K$-dependent entries. This result has an interest in itself.

\subsection{A universality result for symmetric matrices with $K$-dependent
entries}

\label{sectionKdependent}

In this section, we are interested by a universality scheme for the spectral
limiting distribution of symmetric matrices ${\mathbf{X}}_{n}=[X_{k,\ell
}^{(n)}]_{k,\ell=1}^{n}$ normalized by $\sqrt n$ when the entries are
real-valued random variables defined on a common probability space and satisfy
a $K$-dependence condition (see Assumption $\mathbf{A_{3}}$). As we shall see
later, Theorem \ref{thmmdependent} below will be a key step to prove Theorem
\ref{thmUni} .

\smallskip

Let us start by introducing some assumptions concerning the entries
$(X^{(n)}_{k,\ell}, 1 \leq\ell\leq k \leq n)$.

\begin{enumerate}
\item[$\mathbf{A_{1}}$] For all positive integers $n$, ${\mathbb{E}}
(X^{(n)}_{k,\ell}) = 0$ for all $1\leq\ell\leq k \leq n$, and
\[
\frac{1}{n^{2}} \sum_{k=1}^{n}\sum_{\ell=1}^{k} {\mathbb{E}}( | X^{(n)}_{k,
\ell}|^{2} )\leq C < \infty\, .
\]

\item[$\mathbf{A_{2}}$] For any $\tau>0$,
\[
L_{n} ( \tau) := \frac{1}{n^{2}} \sum_{k=1}^{n}\sum_{\ell=1}^{k} {\mathbb{E}}
(| X^{(n)}_{k, \ell}|^{2} \mathbf{1}_{|X^{(n)}_{k, \ell}| > \tau\sqrt{n}})
\rightarrow_{n \rightarrow\infty} 0 \, .
\]

\item[$\mathbf{A_{3}}$] There exists a positive integer $K$ such that for all
positive integers $n$, the following holds: for all nonempty subsets
\[
A,B\subset\{(k,\ell)\in\{1,\dots,n\}^{2}\,|\,1\leq\ell\leq k\leq n\}
\]
such that
\[
\min_{(i,j)\in A}\min_{(k,\ell)\in B}\max(|i-k|,|j-\ell|)>K
\]
the $\sigma$-fields
\[
\sigma\big (X^{(n)}_{i,j}\,,\,(i,j)\in A\big )\ \text{ and }\ \sigma
\big (X^{(n)}_{k,\ell}\,,\,(k,\ell)\in B\big )
\]
are independent.
\end{enumerate}

\bigskip

Condition $\mathbf{A}_{3}$ states that variables with index sets which are at
a distance larger than $K$ are independent.

\bigskip

In Theorem \ref{thmmdependent} below, we then obtain a universality result for
symmetric matrices whose entries are $K$-dependent and satisfy $\mathbf{A}%
_{1}$ and the traditional Lindeberg's condition $\mathbf{A}_{2}$. Note that
$\mathbf{A}_{2}$ is known to be a necessary and sufficient condition for the
empirical spectral distribution of $n^{-1/2}{\mathbf{X}}_{n}$ to converge
almost surely to the semi-circle law when the entries $X_{i,j}^{(n)}$ are
independent, centered and with common variance not depending on $n$ (see
Theorem 9.4.1 in Girko \cite{Girko90}).

\begin{theorem}
\label{thmmdependent} Let ${\mathbf{X}}_{n}=[X_{k,\ell}^{(n)}]_{k,\ell=1}^{n}$
be a symmetric matrix of order $n$ whose entries $(X^{(n)}_{k,\ell}, 1
\leq\ell\leq k \leq n)$ are real-valued random variables satisfying conditions
$\mathbf{A_{1}}$, $\mathbf{A_{2}}$ and $\mathbf{A_{3}}$. Let ${\mathbf{G}}%
_{n}=[G_{i,j}^{(n)}]_{i,j=1}^{n}$ be a symmetric matrix of order $n$ whose
entries $(G_{k,\ell})_{1 \leq\ell\leq k \leq n}$ are real-valued centered
Gaussian random variables with covariance function given by
\begin{equation}
{\mathbb{E}}(G^{(n)}_{k,\ell}G^{(n)}_{i,j})={\mathbb{E}}(X^{(n)}_{k,\ell
}X^{(n)}_{i,j})\,. \label{eqcov1}%
\end{equation}
Then, for any $z\in\mathbb{C}^{+}$,
\begin{equation}
\lim_{n\rightarrow\infty}\big |S_{{\mathbb{X}}_{n}}(z)-{\mathbb{E}%
}(S_{{\mathbb{G}}_{n}}(z)\big )\big |=0\ \text{ almost surely,}
\label{aimthm2}%
\end{equation}
where ${\mathbb{X}}_{n}=n^{-1/2}{\mathbf{X}}_{n}$ and ${\mathbb{G}}%
_{n}=n^{-1/2}{\mathbf{G}}_{n}$.
\end{theorem}
The proof of this result will be given in Appendix.

As we mentioned at the beginning of the section,  this theorem will be a building block to prove that in the stationary and non triangular setting the $K$-dependence condition can be relaxed and more general models for the entries can be considered.  However, the above theorem has also interest in itself. For instance, for the matrices with real
entries, it makes possible to weaken the conditions of Theorems 2.5 and 2.6 in Anderson and Zeitouni
\cite{AZ}. More precisely, due to our Theorem \ref{thmmdependent}, their
assumption 2.2.1 (Ib) can be weaken from the boundness of all moments to the boundness of  moments of order $2$ only 
 plus $\mathbf{A}_{2}$. 
 Furthermore their result
can be strengthened by replacing the convergence in probability to almost sure convergence. Indeed,  our Theorem \ref{thmmdependent} shows that if their
assumption 2.2.1 (Ib) is replaced by $\mathbf{A}_{1}$ plus $\mathbf{A}_{2}$, then to study the limiting spectral distribution we can actually assume without loss of generality that the entries come from a Gaussian random field with the same covariance structure as 
the initial entries. If the $X_{k,\ell}^{(n)}$ are Gaussian random variables then   boundness of all moments means boundness of  moments of order $2$.
\subsection{Proof of Theorem \ref{thmUni}}

For $m$ a positive integer (fixed for the moment) and for any $(u,v)$ in
${\mathbb{Z}}^{2}$ define
\begin{equation}
X_{u,v}^{(m)}={\mathbb{E}}\big (X_{u,v}|{\mathcal{F}}_{u,v}^{(m)}\big )\,,
\label{defXkl}%
\end{equation}
where ${\mathcal{F}}_{u,v}^{(m)}:=\sigma(\xi_{i,j}\,;\,u-m\leq i\leq
u+m,\text{ }v-m\leq j\leq v+m)$.

Let ${\mathbf{X}}_{n}^{(m)}$ be the symmetric random matrix of order $n$
associated with $(X_{u,v}^{(m)})_{(u,v)\in{\mathbb{Z}}^{2}}$ and defined by
$\big ({\mathbf{X}}_{n}^{(m)}\big )_{i,j}=X_{i,j}^{(m)}$ if $1\leq j\leq i\leq
n$ and $\big ({\mathbf{X}}_{n}^{(m)}\big )_{i,j}=X_{j,i}^{(m)}$ if $1\leq
i<j\leq n$. Let
\begin{equation}
{\mathbb{X}}_{n}^{(m)}=n^{-1/2}{\mathbf{X}}_{n}^{(m)}\,. \label{defXklm}%
\end{equation}
We first show that, for any $z\in\mathbb{C}^{+}$,
\begin{equation}
\lim_{m\rightarrow\infty}\limsup_{n\rightarrow\infty}\big |S_{{\mathbb{X}}%
_{n}}(z)-S_{{\mathbb{X}}_{n}^{(m)}}(z)\big |=0\ \text{ a.s.}
\label{approxdependentps}%
\end{equation}
According to Lemma 2.1 in G\"{o}tze \textit{et al.} \cite{GNT} (given for
convenience in Section 5, Lemma \ref{lmagotze}),
\[
\big |S_{{\mathbb{X}}_{n}}(z)-S_{{\mathbb{X}}_{n}^{(m)}}(z)\big |^{2}\leq
\frac{1}{n^{2}v^{4}}\mathrm{Tr}\big (\big ({\mathbf{X}}_{n}-{\mathbf{X}}%
_{n}^{(m)}\big )^{2}\big )\,,
\]
where $v=\operatorname{Im}(z)$. Hence
\[
\big |S_{{\mathbb{X}}_{n}}(z)-S_{{\mathbb{X}}_{n}^{(m)}}(z)\big |^{2}\leq
\frac{2}{n^{2}v^{4}}\sum_{1\leq\ell\leq k\leq n}\big (X_{k,\ell}-X_{k,\ell
}^{(m)}\big )^{2}\,.
\]
Since the shift is ergodic with respect to the measure generated by a sequence
of i.i.d. random variables and the sets of summations are on regular sets, the
ergodic theorem entails that
\[
\lim_{n\rightarrow\infty}\frac{1}{n^{2}}\sum_{1\leq k,\ell\leq n}%
\big (X_{k,\ell}-X_{k,\ell}^{(m)}\big )^{2}={\mathbb{E}}\big (\big (X_{0,0}%
-X_{0,0}^{(m)}\big )^{2}\big )\,\text{ a.s. and in }{\mathbb{L}}^{1}\,.
\]
Therefore%
\begin{equation}
\limsup_{n\rightarrow\infty}\big |S_{{\mathbb{X}}_{n}}(z)-S_{{\mathbb{X}}%
_{n}^{(m)}}(z)\big |^{2}\leq2v^{-4}\Vert X_{0,0}-X_{0,0}^{(m)}\Vert_{2}%
^{2}\,\text{ a.s.} \label{ergothm}%
\end{equation}
Now, by the martingale convergence theorem
\begin{equation}
\Vert X_{0,0}-X_{0,0}^{(m)}\Vert_{2}\rightarrow0\text{ as }m\rightarrow
\infty\,, \label{c1mct}%
\end{equation}
which combined with \eqref{ergothm} proves \eqref{approxdependentps}.

\medskip Let now $(G_{k,\ell}^{(m)})_{(k,\ell)\in{\mathbb{Z}}^{2}}$ be a
real-valued centered Gaussian random field, with covariance function given by
\begin{equation}
{\mathbb{E}}(G_{k,\ell}^{(m)}G_{i,j}^{(m)})={\mathbb{E}}(X_{k,\ell}%
^{(m)}X_{i,j}^{(m)})\text{ for any $(k,\ell)$ and $(i,j)$ in ${\mathbb{Z}}%
^{2}$}\,. \label{egacovfunctionGm}%
\end{equation}
Note that the process $(G_{k,\ell}^{(m)})_{(k,\ell)\in{\mathbb{Z}}^{2}}$ is
then in particular $2m$-dependent. Let now $\mathbf{G}_{n}^{(m)}$ be the
symmetric random matrix of order $n$ defined by $\big (\mathbf{G}_{n}%
^{(m)}\big )_{i,j}=G_{i,j}^{(m)}$ if $1\leq j\leq i\leq n$ and
$\big (\mathbf{G}_{n}^{(m)}\big )_{i,j}=G_{j,i}^{(m)}$ if $1\leq i<j\leq n$.
Denote ${\mathbb{G}}_{n}^{(m)}={\mathbf{G}}_{n}^{(m)}/\sqrt{n}$.

We shall prove that, for any $z\in\mathbb{C}^{+}$,
\begin{equation}
\lim_{n\rightarrow\infty}\big |S_{{\mathbb{X}}_{n}^{(m)}}(z)-{\mathbb{E}%
}\big (S_{{\mathbb{G}}_{n}^{(m)}}(z)\big )\big |=0,\text{ almost surely.}
\label{aimthm2appli}%
\end{equation}
$\ $ With this aim, we shall apply Theorem \ref{thmmdependent} and then show
in what follows that $(X_{k,\ell}^{(m)}\,,\,1\leq\ell\leq k\leq n)$ satisfies
its assumptions.

Note that the sigma-algebras ${\mathcal{F}}_{u,v}^{(m)}:=\sigma(\xi
_{i,j}\,;\,u-m\leq i\leq u+m,\text{ }v-m\leq j\leq v+m)$ and ${\mathcal{F}%
}_{k,\ell}^{(m)}$ are independent as soon as $|u-k|>2m$ or $|v-\ell|>2m$. From
this consideration, we then infer that $(X_{k,\ell}^{(m)}\,,\,1\leq\ell\leq
k\leq n)$ satisfies the assumption $\mathbf{A}_{3}$ of Section
\ref{sectionKdependent} with $K=2m$.

On another hand, since $X_{k,\ell}$ is a centered random variable, so is
$X_{k,\ell}^{(m)}$. Moreover, $\Vert X_{k,\ell}^{(m)}\Vert_{2}\leq\Vert
X_{k,\ell}\Vert_{2}=\Vert X_{1,1}\Vert_{2}$. Hence $(X_{k,\ell}^{(m)}%
\,,\,1\leq\ell\leq k\leq n)$ satisfies the assumption $\mathbf{A}_{1}$ of
Section \ref{sectionKdependent}.

We prove now that the assumption $\mathbf{A}_{2}$ of Section
\ref{sectionKdependent} holds. With this aim, we first notice that, by
Jensen's inequality and stationarity, for any $\tau>0$,
\[
{\mathbb{E}}((X_{k,\ell}^{(m)})^{2}\mathbf{1}_{|X_{k,\ell}^{(m)}|>\tau\sqrt
{n}})\leq{\mathbb{E}}(X_{1,1}^{2}\mathbf{1}_{|X_{1,1}^{(m)}|>\tau\sqrt{n}%
})\,.
\]
Notice now that if $X$ is a real-valued random variable and $\mathcal{F}$ a
sigma-algebra, then for any $\varepsilon>0$,
\[
\mathbb{E}\big (X^{2}\mathbf{1}_{|\mathbb{E}(X|\mathcal{F})|>2\varepsilon
}\big )\leq2\,\mathbb{E}\big (X^{2}\mathbf{1}_{|X|>\varepsilon}\big )\,.
\]
Therefore,
\[
{\mathbb{E}}((X_{k,\ell}^{(m)})^{2}\mathbf{1}_{|X_{k,\ell}^{(m)}|>\tau\sqrt
{n}})\leq2{\mathbb{E}}(X_{1,1}^{2}\mathbf{1}_{|X_{1,1}|>\tau\sqrt{n}/2})\,
\]
which proves that $(X_{k,\ell}^{(m)}\,,\,1\leq\ell\leq k\leq n)$ satisfies
$\mathbf{A}_{2}$ because ${\mathbb{E}}(X_{1,1}^{2})<\infty$.

Since $(X_{k,\ell}^{(m)}\,,\,1\leq\ell\leq k\leq n)$ satisfies the assumptions
$\mathbf{A}_{1}$, $\mathbf{A}_{2}$ and $\mathbf{A}_{3}$ of Section
\ref{sectionKdependent}, applying Theorem \ref{thmmdependent},
\eqref{aimthm2appli} follows.

\medskip According to \eqref{approxdependentps} and \eqref{aimthm2appli}, the
theorem will follow if we prove that, for any $z\in{\mathbb{C}}^{+}$,
\begin{equation}
\lim_{m\rightarrow\infty}\limsup_{n\rightarrow\infty}\big |{\mathbb{E}%
}\big (S_{{\mathbb{G}}_{n}}(z)\big )-{\mathbb{E}}\big (S_{{\mathbb{G}}%
_{n}^{(m)}}(z)\big )\big |=0\,. \label{aim}%
\end{equation}
With this aim, we apply Lemma \ref{interGaussian} from Section
\ref{sectiontechres} which gives
\begin{gather*}
{\mathbb{E}}\big (S_{{\mathbb{G}}_{n}}(z)\big )-{\mathbb{E}}%
\big (S_{{\mathbb{G}}_{n}^{(m)}}(z)\big )\\
=\frac{1}{2}\sum_{1\leq\ell\leq k\leq n}\sum_{1\leq j\leq i\leq n}\int_{0}%
^{1}\big ({\mathbb{E}}(G_{k,\ell}G_{i,j})-{\mathbb{E}}(G_{k,\ell}^{(m)}%
G_{i,j}^{(m)})\big ){\mathbb{E}}\big (\partial_{k\ell}\partial_{ij}%
f({\mathbf{g}}(t))\big ) \, ,
\end{gather*}
where $f$ is defined in \eqref{deffa} and , for $t\in\lbrack0,1]$,
\[
{\mathbf{g}}(t)=(\sqrt{t}G_{k,\ell}+\sqrt{1-t}G_{k,\ell}^{(m)})_{1\leq\ell\leq
k\leq n}\,.
\]
We shall prove that, for any $t$ in $\lbrack0,1]$,
\begin{multline}
\label{estimate}\Big |\sum_{1\leq\ell\leq k\leq n}\sum_{1\leq j\leq i\leq
n}\big ({\mathbb{E}}(G_{k,\ell}G_{i,j})-{\mathbb{E}}(G_{k,\ell}^{(m)}%
G_{i,j}^{(m)})\big ){\mathbb{E}}\big (\partial_{k\ell}\partial_{ij} f
({\mathbf{g}}(t))\big ) \Big |\\
\leq C\Vert X_{0,0}^{(m)}-X_{0,0}\Vert_{2}\Vert X_{0,0}\Vert_{2}\,.
\end{multline}
where $C$ does not depend on $n$ and $t$. After integrating on $[0,1]$ and
then by taking into account that $\Vert X_{0,0}-X_{0,0}^{(m)}\Vert_{2}%
^{2}\rightarrow0$ as $m\rightarrow\infty,$ (\ref{aim}) follows by letting
\ $n$ tend to infinity and then $m$.

To prove (\ref{estimate}), using \eqref{egacovfunctionGm} and
\eqref{egacovfunction}, we write now the following decomposition:
\begin{align}
{\mathbb{E}}(G_{k,\ell}^{{}}G_{i,j}^{{}})-{\mathbb{E}}(G_{k,\ell}^{(m)}%
G_{i,j}^{(m)})  &  ={\mathbb{E}}(X_{k,\ell}^{{}}X_{i,j}^{{}})-{\mathbb{E}%
}(X_{k,\ell}^{(m)}X_{i,j}^{(m)})\nonumber\\
&  = {\mathbb{E}}(X_{k,\ell}( X_{i,j}-X^{(m)}_{i,j}))- {\mathbb{E}}%
((X_{k,\ell}^{(m)}-X_{k,\ell})X_{i,j}^{(m)})\,.\label{approgaussstep2p3}%
\end{align}
We shall decompose the sum on the left-hand side of (\ref{estimate}) in two
sums according to the decomposition (\ref{approgaussstep2p3}) and analyze them
separately. Let us prove that there exists a constant $C$ not depending on $n$
and $t$ such that
\begin{equation}
\Big |\sum_{1\leq\ell\leq k\leq n}\sum_{1\leq j\leq i\leq n}\ {\mathbb{E}%
}((X_{k,\ell}^{(m)}-X_{k,\ell})X_{i,j}^{(m)}){\mathbb{E}}\left(
\partial_{k\ell}\partial_{ij} f ({\mathbf{g}}(t))\right)  \Big | \leq C \Vert
X_{0,0}^{(m)}-X_{0,0}\Vert_{2}\Vert X_{0,0}\Vert_{2}\,.
\label{approgaussstep2p4}%
\end{equation}
To prove \eqref{approgaussstep2p4}, we first notice that without loss of
generality ${\mathbf{g}}(t)$ can be taken independent of $(X_{k,\ell})$ and
then
\[
{\mathbb{E}}((X_{k,\ell}^{(m)}-X_{k,\ell})X_{i,j}^{(m)}){\mathbb{E}}\left(
\partial_{k\ell}\partial_{ij}f({\mathbf{g}}(t))\right)  ={\mathbb{E}}\left(
(X_{k,\ell}^{(m)}-X_{k,\ell})X_{i,j}^{(m)}\partial_{k\ell}\partial_{ij}f
({\mathbf{g}}(t))\,\right)  \, .
\]
Next Lemma \ref{derivatives} from from Section \ref{sectiontechres} applied
with $a_{k,\ell}=(X_{k,\ell}^{(m)}-X_{k,\ell})$ and $b_{k,\ell}=X_{k,\ell
}^{(m)}$ gives: for any $z =u+iv \in{\mathbb{C}}^{+}$,
\begin{align*}
&  \Big |\sum_{1\leq\ell\leq k\leq n}\sum_{1\leq j\leq i\leq n}((X_{k,\ell
}^{(m)}-X_{k,\ell})X_{i,j}^{(m)})\big (\partial_{k\ell}\partial_{ij}f
(\mathbf{g}(t))\big )\Big |\\
&  \leq \frac{2}{v^{3} n^{2}}\Big ( \sum_{1\leq\ell\leq k\leq n}(X_{k,\ell
}^{(m)}-X_{k,\ell})^{2}\Big ) ^{1/2}\Big ( \sum_{1\leq j\leq i\leq n}%
(X_{i,j}^{(m)})^{2}\Big )^{1/2} \, .
\end{align*}
Therefore, by using Cauchy-Schwarz's inequality, we derive
\begin{multline*}
\Big |\sum_{1\leq\ell\leq k\leq n}\sum_{1\leq j\leq i\leq n}\ {\mathbb{E}%
}((X_{k,\ell}^{(m)}-X_{k,\ell})X_{i,j}^{(m)}){\mathbb{E}}\left(
\partial_{k\ell}\partial_{ij}f ({\mathbf{g}}(t))\right)  \Big |\\
\leq\frac{2}{v^{3} n^{2}}\Big (\sum_{1\leq\ell\leq k\leq n}{\mathbb{E}%
}(X_{k,\ell}^{(m)}-X_{k,\ell})^{2}\Big )^{1/2}\Big (\sum_{1\leq j \leq i \leq
n}{\mathbb{E}}(X_{i,j}^{(m)})^{2}\Big )^{1/2}\,.
\end{multline*}
Using stationarity it follows that, for any $z =u+iv \in{\mathbb{C}}^{+}$ and
any $t$ in $\lbrack0,1]$,
\[
\Big |\sum_{1\leq\ell\leq k\leq n}\sum_{1\leq j\leq i\leq n}{\mathbb{E}%
}((X_{k,\ell}^{(m)}-X_{k,\ell})X_{i,j}^{(m)}){\mathbb{E}}\big (\partial
_{k\ell}\partial_{ij} f ({\mathbf{g}}(t))\big )\Big | \leq2 v^{-3} \Vert
X_{0,0}^{(m)}-X_{0,0}\Vert_{2}\Vert X_{0,0}\Vert_{2}\,.
\]
Similarly, we can prove that for any $z =u+iv \in{\mathbb{C}}^{+}$ and any $t$
in $\lbrack0,1]$,
\[
\Big |\sum_{1\leq\ell\leq k\leq n}\sum_{1\leq j\leq i\leq n}{\mathbb{E}%
}(X_{k,\ell}( X_{i,j}-X^{(m)}_{i,j})){\mathbb{E}}\big (\partial_{k\ell
}\partial_{ij} f ({\mathbf{g}}(t))\big )\Big | \leq2 v^{-3} \Vert
X_{0,0}^{(m)}-X_{0,0}\Vert_{2}\Vert X_{0,0}\Vert_{2}\,.
\]
This leads to (\ref{estimate}) and then ends the proof of the theorem.
\hfill$\square$


\subsection{Proof of Theorem \ref{SCL1}}

In order to establish Theorem \ref{SCL1}, it suffices to apply Theorem
\ref{thmUni} and to derive the limit of ${\mathbb{E}}(S_{{\mathbb{G}}_{n}%
}(z))$ for any $z\in{\mathbb{C}}^{+}$, where ${\mathbb{G}}_{n}$ is the
symmetric matrix defined in Theorem \ref{thmUni}. With this aim, we apply
Proposition \ref{propKP} given in Section \ref{sectiontechres}. Proposition
\ref{propKP} is a modification of Theorem 2 in Khorunzhy and Pastur \cite{KP}
(see also Theorem 17.2.1 in \cite{PS}) since in our case, we cannot use
directly the conclusion of their theorem: we are not exactly in the situation
described there. Their symmetric matrix is defined via a \textit{symmetric}
real-valued centered Gaussian random field $(W_{k,\ell})_{k,\ell}$ satisfying
the following property: $W_{k,\ell}=W_{\ell,k}$ for any $(k,\ell
)\in{\mathbb{Z}}^{2}$ and also (2.8) in \cite{KP}. In our situation, and if
\eqref{condgammaY} is assumed, the entries $(g_{k,\ell})_{1\leq k,\ell\leq n}$
of $n^{1/2}{\mathbb{G}}_{n}$ have the following covariances
\begin{equation}
{\mathbb{E}}(g_{i,j}g_{k,\ell})=\gamma_{i-k,j-\ell}(\mathbf{1}_{i\geq
j,k\geq\ell}+\mathbf{1}_{j>i,\ell>k})+\gamma_{i-\ell,j-k}(\mathbf{1}_{i\geq
j,\ell>k}+\mathbf{1}_{j>i,k\geq\ell})\,, \label{equationKP}%
\end{equation}
since by \eqref{egacovfunction} and stationarity
\[
g_{k,\ell}=G_{\max(k,\ell),\min(k,\ell)}\text{ \ and \ }{\mathbb{E}}%
(G_{i,j},G_{k,\ell})=\gamma_{k-i,\ell-j}\,.
\]
Hence, because of the indicator functions appearing in \eqref{equationKP}, our
covariances do not satisfy the condition (2.8) in \cite{KP}. However, the
conclusion of Theorem 2 in \cite{KP} also holds for $S_{{\mathbb{G}}_{n}}(z)$
provided that \eqref{condabssumcov} and \eqref{condgammaY} are satisfied. We
did not find any reference where the assertion above is mentioned so
Proposition \ref{propKP} is proved with this aim. \hfill$\square$

\subsection{Proof of Theorem \ref{thmainGram}}

Let $n=N+p$ and $\mathbb{X}_{n}$ the symmetric matrix of order $n$ defined by
\[
\mathbb{X}_{n}=\frac{1}{\sqrt{p}}\left(
\begin{array}
[c]{cc}%
\mathbf{0}_{p,p} & {\mathcal{X}}_{N,p}^{T}\\
{\mathcal{X}}_{N,p} & \mathbf{0}_{N,N}%
\end{array}
\right)  \,.
\]
Notice that the eigenvalues of $\mathbb{X}^{2}_{n}$ are the eigenvalues of
$p^{-1}{\mathcal{X}}^{T}_{N,p}{\mathcal{X}}_{N,p}$ together with the
eigenvalues of $p^{-1}{\mathcal{X}}_{N,p}{\mathcal{X}}^{T}_{N,p}$. Since these
two latter matrices have the same nonzero eigenvalues, the following relation
holds: for any $z\in{\mathbb{C}}^{+}$
\begin{equation}
S_{\mathbb{B}_{N}}(z)=z^{-1/2}\frac{n}{2N}S_{\mathbb{X}_{n}}(z^{1/2}%
)+\frac{p-N}{2Nz}\,. \label{rel1betweenS}%
\end{equation}
(See, for instance, page 549 in Rashidi Far \textit{et al} \cite{ROBS} for
additional arguments leading to the relation above). A similar relation holds
for the Gram matrix ${\mathbb{H}}_{N}$ associated with the centered Gaussian
random field $(G_{k,\ell})_{(k,\ell)\in{\mathbb{Z}}^{2}}$ having the same
covariance structure as {$(X_{k,\ell})_{(k,\ell)\in{\mathbb{Z}}^{2}}$, namely:
for any $z\in{\mathbb{C}}^{+}$
\begin{equation}
S_{\mathbb{H}_{N}}(z)=z^{-1/2}\frac{n}{2N}S_{\mathbb{G}_{n}}(z^{1/2}%
)+\frac{p-N}{2Nz}\, , \label{rel2betweenS}%
\end{equation}
where $\mathbb{G}_{n}$ is defined as $\mathbb{X}_{n}$ but with ${\mathcal{G}%
}_{N,p}$ replacing ${\mathcal{X}}_{N,p}$. }

In view of the relations \eqref{rel1betweenS} and \eqref{rel2betweenS}, and
since $n/N\rightarrow1+c^{-1}$, to prove \eqref{thmainGramcomp}, it suffices
to show that, for any $z\in{\mathbb{C}}^{+}$,
\begin{equation}
\lim_{n\rightarrow\infty}\big |S_{{\mathbb{X}}_{n}}(z)-{\mathbb{E}%
}\big (S_{{\mathbb{G}}_{n}}(z)\big )\big |=0\ \text{ almost surely}
\label{aim2MPthm}%
\end{equation}
Note now that $\mathbb{X}_{n}:=n^{-1/2}[x_{ij}^{(n)}]_{i,j=1}^{n}$ where
$x_{ij}^{(n)}=\sqrt{\frac{n}{p}}X_{i-p,j}{\mathbf{1}}_{i\geq p+1}{\mathbf{1}%
}_{1\leq j\leq p}$ if $1\leq j\leq i\leq n$, and $x_{ij}^{(n)}=x_{ji}^{(n)}$
if $1\leq i<j\leq n$. Similarly we can write $\mathbb{G}_{n}:=n^{-1/2}%
[g_{ij}^{(n)}]_{i,j=1}^{n}$ where the $g_{ij}^{(n)}$'s are defined as the
$x_{ij}^{(n)}$'s but with $X_{i-p,j}$ replaced by $G_{i-p,j}$. Following the
proof of Theorem \ref{thmUni} , we infer that its conclusion holds (and
therefore \eqref{aim2MPthm} does) even when the stationarity of entries of
$\mathbb{X}_{n}$ and $\mathbb{G}_{n}$ is slightly relaxed as above.
\hfill$\square$

\subsection{Proof of Theorem \ref{MPthm}}

In view of the convergence \eqref{thmainGramcomp}, it suffices to show that
when $N,p\rightarrow\infty$ such that $N/p\rightarrow c \in(0,\infty)$, then
for any $z\in{\mathbb{C}}^{+}$, ${\mathbb{E}} \big ( S_{{\mathbb{H}}_{N}%
}(z)\big ) $ converges to $S(z) =\int_{0}^{1}h(x,z)dx$ where $h(x,z)$ is a
solution to the equation \eqref{eqg}. This follows by applying Theorem 2.1 in
Boutet de Monvel \textit{et al} \cite{BKV} . Indeed setting ${\mathbb{\tilde
H}}_{N} = \frac{p}{N}{\mathbb{H}}_{N}$, this theorem asserts that if
\eqref{condabssumcov} holds then, when $N,p\rightarrow\infty$ such that
$N/p\rightarrow c \in(0,\infty)$, ${\mathbb{E}} \big ( S_{ {\mathbb{\tilde H}%
}_{N}}(z)\big ) $ converges to $m(z)= \int_{0}^{1}v(x,z)dx$, for any
$z\in{\mathbb{C}}^{+}$, where $v(x,z)$ is a solution to the equation
\[
v(x,z)=\Big (-z+c^{-1}\int_{0}^{1}\frac{f(x,s)}{1+ \int_{0}^{1}f(u,s)v(u,z)du}%
ds\Big )^{-1}\, .
\]
This implies that ${\mathbb{E}} \big ( S_{ {\mathbb{H}}_{N}}(z)\big ) $
converges to $S(z)$ as defined in the theorem since the following relation
holds: $S(z) = c^{-1}m(z/c)$.

\hfill$\square$

\section{Technical results}

\label{Sectappendix}

\subsection{Concentration of the spectral measure}

Next proposition is a generalization to row-wise $K$-dependent random matrices
of Theorem 1 (ii) of Guntuboyina and Leeb \cite{GL}.

\begin{proposition}
\label{thconc} Let $(X_{k,\ell}^{(n)})_{1\leq\ell\leq k\leq n}$ be an array of
complex-valued random variables defined on a common probability space. Assume
that there exists a positive integer $K$ such that for any integer
$u\in\lbrack1,n-K]$, the $\sigma$-fields
\[
\sigma\big (X_{i,j}^{(n)}\,,\,1\leq j\leq i\leq u\big )\ \text{ and }%
\ \sigma\big (X_{k,\ell}^{(n)}\,,\,\,1\leq\ell\leq k\,,\,u+K+1\leq k\leq
n\big )
\]
are independent. Define the symmetric matrix ${\mathbb{X}}_{n}$ by
\eqref{defbbX}. Then for every measurable function $f:{\mathbb{R}}%
\rightarrow{\mathbb{R}}$ of bounded variation, any $n\geq K$ and any $r\geq
0$,
\begin{equation}
{\mathbb{P}}\Big (\Big |\int fd\nu_{{\mathbb{X}}_{n}}-{\mathbb{E}}\int
fd\nu_{{\mathbb{X}}_{n}}\Big |\geq r\Big )\leq2\exp\Big (-\frac{nr^{2}%
}{160KV_{f}^{2}}\Big )\,, \label{concmeasure}%
\end{equation}
where $V_{f}$ is the variation of the function $f$.
\end{proposition}

\noindent\textbf{Application to the Stieltjes transform.} Assume that the
assumptions of Proposition \ref{thconc} hold. Let $z=u+{\mathrm{i}}%
v\in{\mathbb{C}}^{+}$ and note that
\[
S_{{\mathbb{X}}_{n}}(z)=\int\frac{1}{x-z}d\nu_{{\mathbb{X}}_{n}}(x)=\int
f_{1}(x)d\nu_{{\mathbb{X}}_{n}}(x)+{\mathrm{i}}\int f_{2}(x)d\nu_{{\mathbb{X}%
}_{n}}(x)\,,
\]
where $f_{1}(x)=\frac{x-u}{(x-u)^{2}+v^{2}}$ and $f_{2}(x)=\frac{v}%
{(x-u)^{2}+v^{2}}$. Now
\[
V_{f_{1}}=\Vert f_{1}^{\prime}\Vert_{1}=\frac{2}{v}\,\text{ and }\,V_{f_{2}%
}=\Vert f_{2}^{\prime}\Vert_{1}=\frac{2}{v}\,.
\]
Therefore, by applying Proposition \ref{thconc} to $f_{1}$ and $f_{2}$, we get
that for any $n\geq K$ and any $r\geq0$,
\begin{equation}
\mathbb{P}\big (\big |S_{{\mathbb{X}}_{n}}(z)-{\mathbb{E}}S_{{\mathbb{X}}_{n}%
}(z)\big |\geq r\big )\leq4\exp\Big (-\frac{nr^{2}v^{2}}{2560K}\Big )\,.
\label{concstieljes}%
\end{equation}

\medskip

\noindent\textbf{Proof of Proposition \ref{thconc}.} It is convenient to start
by considering the map $A$ which "constructs" symmetric matrices of order $n$
as in \eqref{defbbX}. To define it, let $N=n(n+1)/2$ and write elements of
${\mathbb{R}}^{N}$ as $\mathbf{x}=(r_{1},\dots,r_{n})$ where $r_{i}%
=(x_{i,j})_{1\leq j\leq i}$. For any $\mathbf{x}\in{\mathbb{R}}^{N}$, let
$A(\mathbf{x})=A(r_{1},\dots,r_{n})$ be the matrix defined by
\begin{equation}
\big (A(\mathbf{x})\big )_{ij}=\left\{
\begin{array}
[c]{ll}%
\frac{1}{\sqrt{n}}x_{i,j}=\frac{1}{\sqrt{n}}(r_{i})_{j} & \mbox{ if }i\geq j\\
\frac{1}{\sqrt{n}}x_{j,i}=\frac{1}{\sqrt{n}}(r_{j})_{i} & \mbox{ if }i<j
\end{array}
\right.  \label{defofA}%
\end{equation}
For $1\leq i\leq n$, let $R_{i}=(X_{i,j}^{(n)})_{1\leq j\leq i}$. By
definition, we have that ${\mathbb{X}}_{n}=A(R_{1},\dots,R_{n})$. Let then $h$
be the function from ${\mathbb{C}}^{N}$ to ${\mathbb{R}}$ defined by
\[
h(R_{1},\dots,R_{n})=\int fd\nu_{A(R_{1},\dots,R_{n})}\,.
\]
Let $n\geq K$. Denoting by ${\mathcal{F}}_{k}=\sigma(R_{1},\dots,R_{k})$ for
$k\geq1$, and by ${\mathcal{F}}_{0}=\{\emptyset,\Omega\}$, we then write the
following martingale decomposition:
\begin{align*}
\int fd\nu_{{\mathbb{X}}_{n}}  &  -{\mathbb{E}}\int fd\nu_{{\mathbb{X}}_{n}%
}=h(R_{1},\dots,R_{n})-{\mathbb{E}}h(R_{1},\dots,R_{n})\\
&  =\sum_{i=1}^{[n/K]}\Big ({\mathbb{E}}\big (h(R_{1},\dots,R_{n}%
)|{\mathcal{F}}_{iK}\big )-{\mathbb{E}}\big (h(R_{1},\dots,R_{n}%
)|{\mathcal{F}}_{(i-1)K}\big )\Big )\\
&  \quad\quad\quad+{\mathbb{E}}\big (h(R_{1},\dots,R_{n})|{\mathcal{F}}%
_{n}\big )-{\mathbb{E}}\big (h(R_{1},\dots,R_{n})|{\mathcal{F}}_{K[n/K]}%
\big )\\
&  :=\sum_{i=1}^{[n/K]+1}d_{i,n}\,.
\end{align*}
Let
\[
\mathbf{R}_{n}=(R_{1},\dots,R_{n})\ \text{ and }\ \mathbf{R}_{n}^{k,\ell
}=(R_{1},\dots,R_{k},0,\dots,0,R_{\ell+1},\dots,R_{n})\,.
\]
Note now that, for any $i\in\{1,\dots,[n/K]\}$,
\begin{equation}
{\mathbb{E}}\Big (h\big (\mathbf{R}_{n}^{(i-1)K,(i+1)K}\big )|{\mathcal{F}%
}_{iK}\Big )={\mathbb{E}}\Big (h\big (\mathbf{R}_{n}^{(i-1)K,(i+1)K}%
\big )|{\mathcal{F}}_{(i-1)K}\Big )\,. \label{equespcond}%
\end{equation}
To see this it suffices to apply Lemma \ref{lemmaconditional} with
$X=(R_{(i+1)K+1},\dots,R_{n})$, $Y=(R_{1},\dots,R_{(i-1)K})$ and
$Z=(R_{1},\dots,R_{iK})$. Therefore, by taking into account
\eqref{equespcond}, we get that, for any $i\in\{1,\dots,[n/K]\}$,
\begin{multline}
{\mathbb{E}}\big (h(\mathbf{R}_{n})|{\mathcal{F}}_{iK}\big )-{\mathbb{E}%
}\big (h(\mathbf{R}_{n})|{\mathcal{F}}_{(i-1)K}\big )\\
={\mathbb{E}}\big (h(\mathbf{R}_{n})-h\big (\mathbf{R}_{n}^{(i-1)K,(i+1)K}%
\big )|{\mathcal{F}}_{iK}\big )-{\mathbb{E}}\big (h(\mathbf{R}_{n}%
)-h\big (\mathbf{R}_{n}^{(i-1)K,(i+1)K}\big )|{\mathcal{F}}_{(i-1)K}\big )\,.
\label{equespcondbis}%
\end{multline}
We write now that
\begin{multline}
h(\mathbf{R}_{n})-h\big (\mathbf{R}_{n}^{(i-1)K,(i+1)K}%
\big )\label{dec1concen}\\
=\sum_{j=iK+1}^{(i+1)K}\Big (h\big (\mathbf{R}_{n}^{iK,j-1}%
\big )-h\big (\mathbf{R}_{n}^{iK,j}\big )\Big )+\sum_{j=(i-1)K+1}%
^{iK}\Big (h\big (\mathbf{R}_{n}^{j,(i+1)K}\big )-h\big (\mathbf{R}%
_{n}^{j-1,(i+1)K}\big )\Big )\,,
\end{multline}
since $\mathbf{R}_{n}=\mathbf{R}_{n}^{iK,iK}$. But if ${\mathbb{Y}}_{n}$ and
${\mathbb{Z}}_{n}$ are two symmetric matrices of size $n$, then
\[
\Big |\int fd\nu_{{\mathbb{Y}}_{n}}-\int fd\nu_{{\mathbb{Z}}_{n}}\Big |\leq
V_{f}\Vert F^{{\mathbb{Y}}_{n}}-F^{{\mathbb{Z}}_{n}}\Vert_{\infty}%
\]
(see for instance the proof of Theorem 6 in \cite{GL}). Hence, from Theorem
A.43 in Bai and Silverstein \cite{BS},
\[
\Big |\int fd\nu_{{\mathbb{Y}}_{n}}-\int fd\nu_{{\mathbb{Z}}_{n}}%
\Big |\leq\frac{V_{f}}{n}\,\mathrm{rank}\big ({\mathbb{Y}}_{n}-{\mathbb{Z}%
}_{n}\big )\,.
\]
With our notations, this last inequality implies that for any $0\leq k\leq
\ell\leq n$ and $0\leq i\leq j\leq n$
\begin{equation}
\big |h\big (\mathbf{R}_{n}^{k,\ell}\big )-h\big (\mathbf{R}_{n}%
^{i,j}\big )\big |\leq\frac{V_{f}}{n}\,\mathrm{rank}\big (A(\mathbf{R}%
_{n}^{k,\ell})-A(\mathbf{R}_{n}^{i,j})\big )\,. \label{firstinequality}%
\end{equation}
Starting from \eqref{dec1concen} and using \eqref{firstinequality} together
with
\[
\mathrm{rank}\big (A(\mathbf{R}_{n}^{iK,j-1})-A(\mathbf{R}_{n}^{iK,j}%
)\big )\leq2\,
\]
and
\[
\mathrm{rank}\big (A(\mathbf{R}_{n}^{j,(i+1)K})-A(\mathbf{R}_{n}%
^{j-1,(i+1)K})\big )\leq2\,,
\]
we get that
\begin{equation}
\big |h(\mathbf{R}_{n})-h\big (\mathbf{R}_{n}^{(i-1)K,(i+1)K}\big )\big |\leq
\frac{4K}{n}V_{f}\,. \label{ineconc1}%
\end{equation}
Starting from \eqref{equespcondbis} and using \eqref{ineconc1}, it follows
that, for any $i\in\{1,\dots,[n/K]\}$,
\[
\big |{\mathbb{E}}\big (h(\mathbf{R}_{n})|{\mathcal{F}}_{iK}\big )-{\mathbb{E}%
}\big (h(\mathbf{R}_{n})|{\mathcal{F}}_{(i-1)K}\big )\big |\leq\frac{8K}%
{n}V_{f}\newline%
\]
On another hand, since $\mathbf{R}_{n}^{K[n/K],n}$ is ${\mathcal{F}}_{K[n/K]}%
$-measurable,
\begin{multline*}
{\mathbb{E}}\big (h(\mathbf{R}_{n})|{\mathcal{F}}_{n}\big )-{\mathbb{E}%
}\big (h(\mathbf{R}_{n})|{\mathcal{F}}_{K[n/K]}\big )\ ={\mathbb{E}%
}\big (h(\mathbf{R}_{n})-h\big (\mathbf{R}_{n}^{K[n/K],n}\big )|{\mathcal{F}%
}_{n}\big )\\
-{\mathbb{E}}\big (h(\mathbf{R}_{n})-h\big (\mathbf{R}_{n}^{K[n/K],n}%
\big )|{\mathcal{F}}_{K[n/K]}\big )\,.
\end{multline*}
Now
\[
h(\mathbf{R}_{n})-h\big (\mathbf{R}_{n}^{K[n/K],n}\big )=\sum_{j=K[n/K]+1}%
^{n}\Big (h\big (\mathbf{R}_{n}^{j,n}\big )-h\big (\mathbf{R}_{n}%
^{j-1,n}\big )\Big )\,.
\]
So, proceeding as before, we infer that
\[
\big |{\mathbb{E}}\big (h(\mathbf{R}_{n})|{\mathcal{F}}_{n}\big )-{\mathbb{E}%
}\big (h(\mathbf{R}_{n})|{\mathcal{F}}_{K[n/K]}\big )\big |\leq\frac{4K}%
{n}V_{f}\,.
\]
So, overall we derive that $\Vert d_{i,n}\Vert_{\infty}\leq\frac{8K}{n}V_{f}$
for any $i\in\{1,\dots,[n/K]\}$ and $\Vert d_{[n/K]+1,n}\Vert_{\infty}%
\leq\frac{4K}{n}V_{f}$. Therefore, the proposition follows by applying the
Azuma-Hoeffding inequality for martingales. \hfill$\square$

\subsection{Other useful technical results}

\label{sectiontechres}

\begin{lemma}
\label{lemmaconditional}If $X$, $Y$, $Z$ are three random vectors defined on a
probability space $(\Omega,\mathcal{K},\mathbb{P}),$ such that $X$ is
independent of $\sigma(Z)$ and $\sigma(Y)\subset\sigma(Z)$. Then, for any
measurable function $g$ such that $\Vert g(X,Y)\Vert_{1}<\infty$,
\begin{equation}
\mathbb{E}(g(X,Y)|Z)=\mathbb{E}(g(X,Y)|Y)\text{ a.s.} \label{conditional}%
\end{equation}

\end{lemma}

The following lemma is Lemma 2.1 in G\"{o}tze \textit{et al.} \cite{GNT} and
allows to compare two Stieltjes transforms.

\begin{lemma}
\label{lmagotze} Let ${\mathbf{A}}_{n}$ and ${\mathbf{B}}_{n}$ be two
symmetric $n\times n$ matrices. Then, for any $z \in{\mathbb{C}}
\backslash{\mathbb{R}}$,
\[
\big |S_{ {\mathbb{A}}_{n} }(z)- S_{ {\mathbb{B}}_{n} }(z) \big |^{2}
\leq\frac{1}{n^{2}|\operatorname{Im} (z)|^{4}}\mathrm{Tr} \big (({\mathbf{A}%
}_{n} -{\mathbf{B}}_{n} )^{2} \big ) \, ,
\]
where ${\mathbb{A}}_{n} =n^{-1/2} {\mathbf{A}}_{n}$ and ${\mathbb{B}}_{n}
=n^{-1/2} {\mathbf{B}}_{n}$.
\end{lemma}

All along the proofs, we shall use the fact that the Stieltjes transform of
the spectral measure is a smooth function of the matrix entries. Let
$N=n(n+1)/2$ and write elements of ${\mathbb{R}}^{N}$ as $\mathbf{x}%
=(x_{ij})_{1\leq j\leq i\leq n}$. For any $z \in{\mathbb{C}}^{+}$, let
$f(\cdot):=f_{n,z}(\cdot)$ be the function defined from $\mathbb{R}^{N}$ to
$\mathbb{C}$ by
\begin{equation}
f(\mathbf{x})=\frac{1}{n}\mathrm{Tr}\big (A(\mathbf{x})-z{\mathbf{{I}}_{n}%
}\big )^{-1}\ \text{ for any $\mathbf{x}\in\mathbb{R}^{N}$}\,, \label{deffa}%
\end{equation}
where $A(\mathbf{x})$ is the matrix defined in \eqref{defofA} and
$\mathbf{{I}}_{n}$ is the identity matrix of order $n$. The function $f$
admits partial derivatives of all orders. In particular, denoting for any
$\mathbf{u}\in\{(i,j)\}_{1\leq j\leq i\leq n}$, $\partial_{\mathbf{u}}f$ for
$\partial f_{n}/\partial x_{\mathbf{u}}$, the following upper bounds hold: for
any $z=x+iy\in\mathbb{C}^{+}$ and any $\mathbf{u},\mathbf{v},\mathbf{w}$ in
$\{(i,j)\}_{1\leq j\leq i\leq n}$,
\begin{equation}
|\partial_{\mathbf{u}}f|\leq\frac{2}{y^{2} n^{3/2}}\,,\,|\partial_{\mathbf{u}%
}\partial_{\mathbf{v}}f|\leq\frac{4}{y^{3} n^{2}}\,\text{ and }\,|\partial
_{\mathbf{u}}\partial_{\mathbf{v}}\partial_{\mathbf{w}} |\leq\frac
{3\times2^{5/2}}{y^{4} n^{5/2}}\,. \label{boundsd}%
\end{equation}
(See the equalities (20) and (21) in \cite{Ct} together with the computations
on pages 2074-2075). In addition, the following lemma has been proved in
Merlev\`ede and Peligrad (2014):

\begin{lemma}
\label{derivatives} Let $z=x+iy\in\mathbb{C}^{+}$ and $f_{n}:=f_{n,z}$ be
defined by \eqref{deffa}. Let $(a_{ij})_{1\leq j\leq i\leq n}$ and
$(b_{ij})_{1\leq j\leq i\leq n}$ be real numbers. Then, for any subset
${\mathcal{I}}_{n}$ of $\{(i,j)\}_{1\leq j\leq i\leq n}$ and any element
${\mathbf{x}}$ of ${\mathbb{R}}^{N}$,
\[
\Big |\sum_{\mathbf{u}\in{\mathcal{I}}_{n}}\sum_{\mathbf{v}\in{\mathcal{I}%
}_{n}}a_{\mathbf{u}}b_{\mathbf{v}}\partial_{\mathbf{u}}\partial_{\mathbf{v}%
}f({\mathbf{x}})\Big |\leq\frac{2}{y^{3} n^{2}}\Big (\sum_{\mathbf{u}%
\in{\mathcal{I}}_{n}}a_{\mathbf{u}}^{2}\sum_{\mathbf{v}\in{\mathcal{I}}_{n}%
}b_{\mathbf{v}}^{2}\Big )^{1/2}\,.
\]

\end{lemma}

Next lemma is a consequence of the well-known Gaussian interpolation trick.

\begin{lemma}
\label{interGaussian}Let $(Y_{k,\ell})_{(k,\ell)\in{\mathbb{Z}}^{2}}$ and
$(Z_{k,\ell})_{(k,\ell)\in{\mathbb{Z}}^{2}}$ be two centered real-valued
Gaussian processes. Let ${\mathbf{Y}}_{n}$ be the symmetric random matrix of
order $n$ defined by $\big ({\mathbf{Y}}_{n}\big )_{i,j}=Y_{i,j}$ if $1\leq
j\leq i\leq n$ and $\big ({\mathbf{Y}}_{n}\big )_{i,j}=Y_{j,i}$ if $1\leq
i<j\leq n$. Denote ${\mathbb{Y}}_{n}=\frac{1}{\sqrt{n}}{\mathbf{Y}}_{n}$.
Define similarly ${\mathbb{Z}}_{n}$. Then, for any $z = x+{\mathrm{i}}y
\in{\mathbb{C}^{+}}$,
\begin{align}
\label{res1lmainterGaussian}{\mathbb{E}} \big ( S_{{\mathbb{Y}}_{n}}(z)
\big )  &  -{\mathbb{E}}\big (S_{{\mathbb{Z}}_{n}}(z)\big )\nonumber\\
&  = \frac{1}{2}\sum_{1\leq\ell\leq k\leq n}\sum_{1\leq j\leq i\leq n}
\int_{0}^{1} \big ( {\mathbb{E}} (Y_{k,\ell}Y_{i,j}) - {\mathbb{E}}
(Z_{k,\ell}Z_{i,j})\big ) {\mathbb{E}} \big ( \partial_{k \ell} \partial_{ij}
f ( {\mathbf{u}} (t) ) \big ) dt
\end{align}
where, for $t \in[0,1]$, ${\mathbf{u}} (t)= (\sqrt{t} Y_{k,\ell} + \sqrt{1-t}
Z_{k,\ell})_{1 \leq\ell\leq k \leq n}$ and
\begin{equation}
\label{res2lmainterGaussian}\big |{\mathbb{E}} \big ( S_{{\mathbb{Y}}_{n}}(z)
\big )-{\mathbb{E}}\big (S_{{\mathbb{Z}}_{n}}(z)\big )\big |\leq\frac{2}{n^{2}
y^{3}}\sum_{1\leq\ell\leq k\leq n}\sum_{1\leq j\leq i\leq n}| {\mathbb{E}}
(Y_{k,\ell}Y_{i,j}) - {\mathbb{E}} (Z_{k,\ell}Z_{i,j})|\,.
\end{equation}

\end{lemma}

\noindent\textbf{Proof.} Using the definition of $f$, we first write
\[
{\mathbb{E}} \big ( S_{{\mathbb{Y}}_{n}}(z) \big ) = {\mathbb{E}} f
\big ( (Y_{k,\ell})_{1 \leq\ell\leq k \leq n} \big )\text{ and }
\ {\mathbb{E}}\big (S_{{\mathbb{Z}}_{n}}(z)\big ) = {\mathbb{E}} f
\big ( (Z_{k,\ell})_{1 \leq\ell\leq k \leq n} \big ) \, .
\]
Equality \eqref{res1lmainterGaussian} then follows from the usual
interpolation trick (for an easy reference we cite Talagrand \cite{Tal}
Section 1.3, Lemma 1.3.1.). To obtain the upper bound
\eqref{res2lmainterGaussian}, it suffices then to take into account
\eqref{boundsd}. \hfill$\square$

\medskip

Below we give a Taylor expansion for functions of random variables of a
convenient type for Lindeberg's method.

\begin{lemma}
\label{Taylor variables}Let $f(\cdot)$ be a function from ${\mathbb{R}}^{d+m}$
to ${\mathbb{C}}$, three times differentiable, with continuous and bounded
third partial derivatives, i.e. there is a constant $L_{3}$ such that
\[
\vert\partial_{i}\partial_{j}\partial_{k}f ({\mathbf{x}}) \vert\leq
L_{3}\text{ for all $i,j,k$ and ${\mathbf{x}}$} \, .
\]
Let $\mathbf{X},\mathbf{Y}, \mathbf{Z}$ be random vectors defined on a
probability space $(\Omega,\mathcal{K},\mathbb{P})$, such that $\mathbf{X}$
and $\mathbf{Y}$ take their values in $\mathbb{R}^{d}$, and $\mathbf{Z}$ takes
its values in $\mathbb{R}^{m}$. Assume in addition that $\mathbf{X}$ and
$\mathbf{Y}$ are independent of $\mathbf{Z}$, and that $\mathbf{X}$ and
$\mathbf{Y}$ are in ${\mathbb{L}}^{3} ({\mathbb{R}}^{d})$, centered at
expectation and have the same covariance structure. Then, for any permutation
$\pi:{\mathbb{R}}^{d+m}\rightarrow{\mathbb{R}}^{d+m}$, we have
\[
|\mathbb{E}f(\pi(\mathbf{X},\mathbf{Z}))-\mathbb{E}f(\pi(\mathbf{Y}%
,\mathbf{Z}))|\leq\frac{L_{3}d^{2}}{3}\Big (  \sum_{j=1}^{d} {\mathbb{E}}
(|X_{j}|^{3}) +\sum_{j=1}^{d}{\mathbb{E}} (|Y_{j}|^{3})\Big )  \, .
\]

\end{lemma}

The proof of this lemma is based on the following Taylor expansion for
functions of several variables.

\begin{lemma}
\label{Taylor}Let $g(\cdot)$ be a function from ${\mathbb{R}}^{p}$ to
${\mathbb{R}}$, three times differentiable, with continuous third partial
derivatives and such that
\[
|\partial_{i}\partial_{j}\partial_{k}g({\mathbf{x}})|\leq L_{3}\text{ for all
$i,j,k$ and ${\mathbf{x}}$}\,.
\]
Then, for any $\mathbf{a}=(a_{1},\dots,a_{p})$ and $\mathbf{b}=(b_{1}%
,\dots,b_{p})$ in ${\mathbb{R}}^{p}$,
\begin{gather*}
g(\mathbf{b})-g(\mathbf{a})=\sum_{k=1}^{p}(b_{j}-a_{j})\partial_{j}%
g(\mathbf{0})+\frac{1}{2}\sum_{j,k=1}^{p}(b_{j}b_{k}-a_{j}a_{k})\partial
_{j}\partial_{k}g(\mathbf{0})+R_{3}(\mathbf{a},\mathbf{b}).\\
\text{with }|R_{3}(\mathbf{a},\mathbf{b})|\leq\frac{L_{3}}{6}\Big (\big (\sum
_{j=1}^{p}|a_{j}|\big )^{3}+\big (\sum_{j=1}^{p}|b_{j}|\big )^{3}%
\Big )\leq\frac{L_{3}p^{2}}{6}\Big (\sum_{j=1}^{p}|a_{j}|^{3}+|b_{j}%
|^{3}\Big)\,.
\end{gather*}

\end{lemma}

\noindent\textbf{Proof of Lemma \ref{Taylor}.} We use Taylor expansion of
second order for functions with bounded partial derivatives of order three. It
is well-known that
\begin{gather*}
g(\mathbf{a})-g(\mathbf{0}_{p})=\sum_{j=1}^{p}a_{j}\partial_{j}g(\mathbf{0}%
_{p})+\frac{1}{2}\sum_{j,k=1}^{p}a_{j}a_{k}\partial_{j}\partial_{k}%
g(\mathbf{0}_{p})+R_{3}(\mathbf{a})\,,\\
\text{where }|R_{3}(\mathbf{a})|\leq\frac{L_{3}}{6} \big (\sum_{j=1}^{p}%
|a_{j}| \big )^{3}\leq\frac{L_{3}p^{2}}{6}\sum_{j=1}^{p}|a_{j}|^{3}\,.
\end{gather*}
By writing a similar expression for $g(\mathbf{b})-g(\mathbf{0}_{p})$ and
substracting them the result follows. \hfill$\square$

\bigskip

\noindent\textbf{Proof of Lemma \ref{Taylor variables}. }For simplicity of the
notation we shall prove it first for $f((\mathbf{X},\mathbf{Z}))-f((\mathbf{Y}%
,\mathbf{Z})).$ We start by applying Lemma \ref{Taylor} to real and imaginary
part of $f$ and obtain
\begin{gather*}
f(\mathbf{X},\mathbf{Z})-f(\mathbf{Y},\mathbf{Z})=\sum_{j=1}^{d}(X_{j}%
-Y_{j})\partial_{j}f(\mathbf{0}_{d},\mathbf{Z})+\frac{1}{2}\sum_{j,k=1}%
^{d}(X_{k}X_{j}-Y_{k}Y_{j})\partial_{j}\partial_{k}f(\mathbf{0}_{d}%
,\mathbf{Z})+R_{3} \, ,\\
\text{with }|R_{3}|\leq\frac{L_{3}d^{2}}{3}\Big ( \sum_{j=1}^{d}|X_{j}%
|^{3}+\sum_{j=1}^{d}|Y_{j}|^{3}\Big) \, .
\end{gather*}
By taking the expected value and taking into account the hypothesis of
independence and the fact that $\mathbf{X}$ and $\mathbf{Y}$ are centered at
expectations and have the same covariance structure, we obtain, for all $1\leq
j\leq d$
\[
\mathbb{E(}(X_{j}-Y_{j})\partial_{j}f(\mathbf{0}_{d},\mathbf{Z}))=(\mathbb{E}%
X_{j}-\mathbb{E}Y_{j})\mathbb{E}\partial_{j}f(\mathbf{0}_{d},\mathbf{Z})=0
\]
and, for all $1\leq k,j\leq d$,
\[
\mathbb{E}(X_{k}X_{j}-Y_{k}Y_{j})\partial_{j}\partial_{k}f(\mathbf{0}_{d}
,\mathbf{Z})=(\mathbb{E(}X_{k}X_{j})-\mathbb{E(}Y_{k}Y_{j}))\mathbb{E}%
\partial_{j}\partial_{k}f(\mathbf{0}_{d},\mathbf{Z})=0\, .
\]
It follows that
\begin{gather*}
\mathbb{E}f(\mathbf{X},\mathbf{Z})-\mathbb{E}f(\mathbf{Y},\mathbf{Z})=R_{3} \,
,\\
\text{with }|R_{3}|\leq\frac{L_{3}d^{2}}{3}\Big ( \sum_{j=1}^{d}|X_{j}%
|^{3}+\sum_{j=1}^{p}|Y_{j}|^{3}\Big) \, .
\end{gather*}
It remains to note that the result remains valid for any permutation of
variables ($\mathbf{X},\mathbf{Z}$). The variables in $\mathbf{X},\mathbf{Z}$
can hold any positions among the variables in function $f$ since we just need
all the derivatives of order three to be uniformly bounded. The difference in
the proof consists only in re-denoting the partial derivatives; for instance
instead of $\partial_{j}$ we shall use $\partial_{k_{j}}$ where $k_{j}$ ,
$1\leq k_{j}\leq d+m$ denotes$\ $the index of the variable $X_{j}$ in
$f(x_{1},x_{2},...,x_{d+m}).$ \hfill$\square$

\bigskip

We provide next a technical lemma on the behavior of the expected value of
Stieltjes transform of symmetric matrices with Gaussian entries. In Lemma
\ref{lmaKP} and Proposition \ref{propKP} below, we consider a stationary
real-valued centered Gaussian random field $(G_{k,\ell})_{(k,\ell
)\in{\mathbb{Z}}^{2}}$ with covariance function given by: for any $(k,\ell
)\in{\mathbb{Z}}^{2}$ and any $(i,j)\in{\mathbb{Z}}^{2}$,
\[
{\mathbb{E}}(G_{k,\ell} G_{i,j})=\gamma_{k-i,\ell-j}\,,
\]
satisfying \eqref{condabssumcov} and \eqref{condgammaY}. We define then two
symmetric matrices of order $n$, ${\mathbb{G}}_{n}=n^{-1/2}[g_{k,\ell
}]_{k,\ell=1}^{n}$ and ${\mathbb{W}}_{n}=n^{-1/2}[W_{k,\ell}]_{k,\ell=1}^{n}$
where the entries $g_{k,\ell}$ and $W_{k,\ell}$ are defined respectively by
\[
g_{k,\ell}=G_{\max(k,\ell),\min(k,\ell)} \ \text{ and } \ W_{k,\ell} =
\frac{1}{\sqrt{2}} \big ( G_{k,\ell} + G_{\ell,k}\big ) \, .
\]

\begin{lemma}
\label{lmaKP} For any $z \in{\mathbb{C}} \backslash{\mathbb{R}}$ the following
convergence holds:
\[
\lim_{n \rightarrow\infty} \big |{\mathbb{E}} ( S_{{\mathbb{G}}_{n}}(z) ) -
{\mathbb{E}} ( S_{{\mathbb{W}}_{n}}(z) ) \big | = 0 \, .
\]

\end{lemma}

As a consequence of this lemma and Theorem 2 in \cite{KP}, we obtain the
following result concerning the limiting spectral distribution of both
${\mathbb{G}}_{n}$ and ${\mathbb{W}}_{n}$.

\begin{proposition}
\label{propKP} For any $z \in{\mathbb{C}} \backslash{\mathbb{R}}$,
$S_{{\mathbb{G}}_{n}}(z)$ and $S_{{\mathbb{W}}_{n}}(z)$ have almost surely the
same limit, $S(z)$, defined by the relations \eqref{defSlim} and \eqref{defSlim2}.
\end{proposition}

\noindent\textbf{Proof of Lemma \ref{lmaKP}.} According to Lemma
\ref{interGaussian}, for any $z\in{\mathbb{C}}\backslash{\mathbb{R}}$,
\[
\big |{\mathbb{E}}\big (S_{{\mathbb{G}}_{n}}(z)\big )-{\mathbb{E}%
}\big (S_{{\mathbb{W}}_{n}}(z)\big )\big |\leq\frac{2}{n^{2}|\operatorname{Im}%
(z)|^{3}}\sum_{1\leq\ell\leq k\leq n}\sum_{1\leq j\leq i\leq n}|\mathrm{Cov}%
(G_{k,\ell},G_{i,j})-\mathrm{Cov}(W_{k,\ell},W_{i,j})|\,.
\]
Taking into account \eqref{condgammaY}, we get
\begin{equation}
{\mathbb{E}}(W_{k,\ell}W_{i,j})=\gamma_{k-i,\ell-j}+\gamma_{k-j,\ell
-i}\,.\label{covdesWij}%
\end{equation}
Hence,
\[
\big |{\mathbb{E}}\big (S_{{\mathbb{G}}_{n}}(z)\big )-{\mathbb{E}%
}\big (S_{{\mathbb{W}}_{n}}(z)\big )\big |\leq\frac{2}{n^{2}|\operatorname{Im}%
(z)|^{3}}\sum_{1\leq\ell\leq k\leq n}\sum_{1\leq j\leq i\leq n}|\gamma
_{k-j,\ell-i}|\,.
\]
Using \eqref{condgammaY} and noticing that by stationarity $\gamma
_{u,v}=\gamma_{-u,-v}$ for any $(u,v)\in{\mathbb{Z}}^{2}$, we
get
\[
\sum_{1\leq\ell\leq k\leq n}\sum_{1\leq j\leq i\leq n}|\gamma_{k-j,\ell
-i}|\leq2\sum_{1\leq\ell\leq k\leq n}\sum_{1\leq j\leq i\leq k}|\gamma
_{k-j,\ell-i}|\,.
\]
By simple algebra, we infer that, for any positive integer $m_{n}$ less than
$n$,
\[
\big |{\mathbb{E}}\big (S_{{\mathbb{G}}_{n}}(z)\big )-{\mathbb{E}%
}\big (S_{{\mathbb{W}}_{n}}(z)\big )\big |\leq\frac{4}{|\operatorname{Im}%
(z)|^{3}}\Big (\frac{2m_{n}}{n}\sum_{p=0}^{n}\sum_{q=-n}^{n}|\gamma
_{p,q}|+\sum_{p\geq m_{n}}\sum_{q\in{\mathbb{Z}}}|\gamma_{p,q}|\Big )\,,
\]
for any $z\in{\mathbb{C}}\backslash{\mathbb{R}}$. The lemma then follows by
taking into account \eqref{condabssumcov} and by selecting $m_{n}$ such that
$m_{n}\rightarrow\infty$ and $m_{n}/n\rightarrow0$. \hfill$\square$

\medskip

\noindent\textbf{Proof of Proposition \ref{propKP}.} The Borel-Cantelli lemma
together with Theorem 17.1.1 in \cite{PS} imply that, for any $z
\in{\mathbb{C}} \backslash{\mathbb{R}}$,
\[
\text{$\lim_{n \rightarrow\infty} |S_{{\mathbb{G}}_{n}}(z) - {\mathbb{E}}
\big (S_{{\mathbb{G}}_{n}}(z) \big ) | = 0$ and $\lim_{n \rightarrow\infty}
|S_{{\mathbb{W}}_{n}}(z) - {\mathbb{E}} \big (S_{{\mathbb{W}}_{n}}(z) \big ) |
= 0$ a.s.}
\]
Therefore, the proposition follows by Lemma \ref{lmaKP} combined with Theorem
2 in \cite{KP} applied to ${\mathbb{E}} \big (S_{{\mathbb{W}}_{n}}(z)\big )$.
Indeed the entries $(W_{k,\ell})_{1 \leq k,\ell\leq n}$ of the matrix
$n^{1/2}{\mathbb{W}}_{n}$ form a \textit{symmetric} real-valued centered
Gaussian random field whose covariance function satisfies \eqref{covdesWij}.
Hence relation (2.8) in \cite{KP} holds. In addition, by \eqref{condabssumcov}
, condition (2.9) in \cite{KP} is also satisfied. At this step, the reader
should notice that Theorem 2 in \cite{KP} also requires additional conditions
on the covariance function $\gamma_{k,\ell}$ (this function is denoted by
$B(k,\ell)$ in this latter paper), namely $\gamma_{k,\ell}=\gamma_{\ell
,k}=\gamma_{\ell,-k}$. In our case, the first  holds (this is
\eqref{condgammaY}) but not necessarily $\gamma_{\ell,k}=\gamma_{\ell,-k}$
since by stationarity we only have $\gamma_{\ell,k}=\gamma_{-\ell,-k}$.
However a careful analysis of the proof of Theorem 2 in \cite{KP} (and in
particular of their auxiliary lemmas) or of the proof of Theorem 17.2.1 in
\cite{PS}, shows that the only condition required on the covariance function
to derive the limiting equation of the Stieljes transform is the absolute
summability condition (2.9) in \cite{KP}. It is noteworthy to indicate that,
in Theorem 2 of \cite{KP}, the symmetry conditions on the covariance function
$\gamma_{k,\ell}$ must only translate the fact that the entries of the matrix
form a stationary \textit{symmetric} real-valued centered Gaussian random
field, so $\gamma_{k,\ell}$ has only to satisfy $\gamma_{k,\ell}=\gamma
_{\ell,k}=\gamma_{-\ell,-k}$ for any $(k,\ell) \in{\mathbb{Z}}^{2}$.
\hfill$\square$

\section{Appendix: proof of Theorem \ref{thmmdependent}}

By using inequality \eqref{concstieljes} together with the Borel-Cantelli
lemma, it follows that, for any $z\in\mathbb{C}^{+}$,
\[
\lim_{n\rightarrow\infty}\big |S_{{\mathbb{X}}_{n}}(z)-{\mathbb{E}%
}\big (S_{{\mathbb{X}}_{n}}(z)\big )\big |=0\ \text{ almost surely.}%
\]
To prove the almost sure convergence \eqref{aimthm2}, it suffices then to
prove that, for any $z\in\mathbb{C}^{+}$,
\begin{equation}
\lim_{n\rightarrow\infty}\big |{\mathbb{E}}\big (S_{{\mathbb{X}}_{n}%
}(z)\big )-{\mathbb{E}}\big (S_{{\mathbb{G}}_{n}}(z)\big )\big |=0\,.
\label{aimthm1*}%
\end{equation}
We start by truncating the entries of the matrix ${\mathbb{X}}_{n}$. Since
$\mathbf{A_{2}}$ holds, we can consider a decreasing sequence of positive
numbers $\tau_{n}$ such that, as $n\rightarrow\infty$,
\begin{equation}
\tau_{n}\rightarrow0\ ,\ L_{n}(\tau_{n})\rightarrow0\ \text{ and }\ \tau
_{n}\sqrt{n}\rightarrow\infty\,. \label{choicepntaun}%
\end{equation}
Let ${\bar{\mathbf{X}}}_{n}=[{\bar{X}}^{(n)}_{k,\ell}]_{k,\ell=1}^{n}$ be the
symmetric matrix of order $n$ whose entries are given by:
\[
{\bar{X}}^{(n)}_{k,{\ell}}=X^{(n)}_{k,{\ell}}\mathbf{1}_{|X^{(n)}_{k,\ell
}|\leq\tau_{n}\sqrt{n}}-{\mathbb{E}}\big (X^{(n)}_{k,{\ell}}\mathbf{1}%
_{|X^{(n)}_{k,\ell}|\leq\tau_{n}\sqrt{n}}\big )\, .
\]
Define ${\bar{\mathbb{X}}}_{n}:=n^{-1/2}{\bar{\mathbf{X}}}_{n}$. Using
\eqref{choicepntaun}, it has been proved in Section 2.1 of \cite{GNT} that,
for any $z\in\mathbb{C}^{+}$,
\[
\lim_{n\rightarrow\infty}\big |{\mathbb{E}}\big (S_{{\mathbb{X}}_{n}%
}(z)\big )-{\mathbb{E}}\big (S_{{\bar{\mathbb{X}}}_{n}}(z)\big )\big |=0\,.
\]
Therefore, to prove \eqref{aimthm1*} (and then the theorem), it suffices to
show that, for any $z\in\mathbb{C}^{+}$,
\begin{equation}
\lim_{n\rightarrow\infty}\big |{\mathbb{E}}\big (S_{{\bar{\mathbb{X}}}_{n}%
}(z)\big )-{\mathbb{E}}\big (S_{{\mathbb{G}}_{n}}(z)\big )\big |=0\,.
\label{aimthm1}%
\end{equation}

The proof of \eqref{aimthm1} is then divided in three steps. The first step
consists of replacing in \eqref{aimthm1}, the matrix ${\mathbb{G}}_{n}$ by a
symmetric matrix ${\bar{\mathbb{G}}}_{n}$ of order $n$ whose entries are
real-valued Gaussian random variables with the same covariance structure as
the entries of ${\bar{\mathbb{X}}}_{n}$. The second step consists of
"approximating" ${\bar{\mathbb{X}}}_{n}$ and ${\bar{\mathbb{G}}}_{n}$ by
matrices with "big square independent blocks" containing the entries spaced by
"small blocks" around them containing only zeros as entries. Due to the
assumption $\mathbf{A_{3}}$, the random variables contained in two different
big blocks will be independent. The third and last step consists of proving
the mean convergence \eqref{aimthm1} but with ${\bar{\mathbb{X}}}_{n}$ and
${\mathbb{G}}_{n}$ replaced by their approximating matrices with independent
blocks. This step will be achieved with the help of the Lindeberg method.
\bigskip

\noindent\textit{Step 1.} Let ${\bar{\mathbf{G}}}_{n}=[{\bar{G}}^{(n)}%
_{k,\ell}]_{k,\ell=1}^{n}$ be the symmetric matrix of order $n$ whose entries
$({\bar{G}}^{(n)}_{k,\ell} \, , \, 1 \leq\ell\leq k \leq n)$ are real-valued
centered Gaussian random variables with the following covariance structure:
for any $1 \leq\ell\leq k \leq n$ and any $1 \leq j \leq i \leq n$,
\begin{equation}
{\mathbb{E}}({\bar{G}}^{(n)}_{k,\ell}{\bar{G}}^{(n)}_{i,j})={\mathbb{E}}%
({\bar{X}}^{(n)} _{k,\ell}{\bar{X}}^{(n)}_{i,j}) \,.
\label{egacovfunctionGmtruncated}%
\end{equation}
There is no loss of generality by assuming in the rest of the proof that the
$\sigma$-fields $\sigma( {\bar{G}}^{(n)}_{k,\ell} \, , \, 1 \leq\ell\leq k
\leq n)$ and $\sigma( X^{(n)}_{k,\ell} \, , \, 1 \leq\ell\leq k \leq n)$ are independent.

Denote ${\bar{\mathbb{G}}}_{n} = \frac{1}{\sqrt{ n} }{\bar{\mathbf{G}}}_{n}$.
We shall prove that, for any $z\in\mathbb{C}^{+}$,
\begin{equation}
\lim_{n\rightarrow\infty}\big |{\mathbb{E}}\big (S_{{\bar{\mathbb{G}}}_{n}%
}(z)\big )-{\mathbb{E}}\big (S_{{\mathbb{G}}_{n}}(z)\big )\big |=0\,.
\label{aimthm1approxgau}%
\end{equation}
Applying Lemma \ref{interGaussian}, we get
\begin{equation}
\big |{\mathbb{E}}\big (S_{{\bar{\mathbb{G}}}_{n}}(z)\big )-{\mathbb{E}%
}\big (S_{{\mathbb{G}}_{n}}(z)\big )\big |\leq\frac{2}{v^{3}n^{2}}\sum
_{1\leq\ell\leq k\leq n}\sum_{1\leq j\leq i\leq n}|{\mathbb{E}}(G^{(n)}%
_{k,\ell}G^{(n)}_{i,j})-{\mathbb{E}}({\bar{G}}^{(n)}_{k,\ell}{\bar{G}}%
^{(n)}_{i,j})|\, , \label{approxtrungau}%
\end{equation}
where $v = \operatorname{Im} (z)$. Recall now that ${\mathbb{E}}%
(G^{(n)}_{k,\ell}G^{(n)}_{i,j})={\mathbb{E}}(X^{(n)}_{k,\ell}X^{(n)}_{i,j})$
and ${\mathbb{E}}({\bar{G}}^{(n)} _{k,\ell}{\bar{G}}^{(n)}_{i,j})={\mathbb{E}%
}({\bar{X}}^{(n)}_{k,\ell}{\bar{X}}^{(n)}_{i,j})$. Hence, setting $b_{n}%
=\tau_{n}\sqrt{n}$, we have
\[
{\mathbb{E}}(G^{(n)}_{k,\ell}G^{(n)}_{i,j}) -{\mathbb{E}}({\bar{G}}%
^{(n)}_{k,\ell}{\bar{G}}^{(n)}_{i,j})=\mathrm{Cov}\big (X^{(n)}_{k,\ell
},X^{(n)}_{i,j}\big ) -\mathrm{Cov}\big (X^{(n)}_{k,\ell}\mathbf{1}%
_{|X^{(n)}_{k,\ell}|\leq b_{n}},X^{(n)}_{i,j}\mathbf{1}_{|X^{(n)}_{i,j}|\leq
b_{n}}\big )\,.
\]
Note that
\begin{gather*}
\mathrm{Cov}\big (X^{(n)}_{k,\ell},X^{(n)}_{i,j}\big ) -\mathrm{Cov}%
\big (X^{(n)}_{k,\ell}\mathbf{1}_{|X^{(n)}_{k,\ell}|\leq b_{n}},X^{(n)}%
_{i,j}\mathbf{1}_{|X^{(n)}_{i,j}|\leq b_{n}}\big )\\
=\mathrm{Cov}\big (X^{(n)}_{k,\ell}\mathbf{1}_{|X^{(n)}_{k,\ell}|\leq b_{n}%
},X^{(n)}_{i,j}\mathbf{1}_{|X^{(n)}_{i,j}|>b_{n}}\big )+\mathrm{Cov}%
\big (X^{(n)}_{k,\ell}\mathbf{1}_{|X^{(n)}_{k,\ell}|>b_{n}},X^{(n)}%
_{i,j}\mathbf{1}_{|X^{(n)}_{i,j}|>b_{n}}\big )\\
+\mathrm{Cov}\big (X^{(n)}_{k,\ell}\mathbf{1}_{|X^{(n)}_{k,\ell}|>b_{n}%
},X^{(n)}_{i,j}\mathbf{1}_{|X^{(n)}_{i,j}|\leq b_{n}}\big )\,
\end{gather*}
implying, by Cauchy-Schwarz's inequality, that
\begin{align*}
&  \big |\mathrm{Cov}\big (X^{(n)}_{k,\ell},X^{(n)}_{i,j}\big ) -\mathrm{Cov}%
\big (X^{(n)}_{k,\ell}\mathbf{1}_{|X^{(n)}_{k,\ell}|\leq b_{n}},X^{(n)}%
_{i,j}\mathbf{1}_{|X^{(n)}_{i,j}|\leq b_{n}}\big )\big |\\
&  \leq2b_{n}{\mathbb{E}}\big (|X^{(n)}_{i,j}|\mathbf{1}_{|X^{(n)}%
_{i,j}|>b_{n}}\big )+2b_{n}{\mathbb{E}}\big (|X^{(n)}_{k,\ell}|\mathbf{1}%
_{|X^{(n)}_{k,\ell}|>b_{n}}\big )+2\Vert X^{(n)}_{i,j}\mathbf{1}%
_{|X^{(n)}_{i,j}|>b_{n}}\Vert_{2}\Vert X^{(n)}_{k,\ell}\mathbf{1}%
_{|X^{(n)}_{k,\ell}|>b_{n}}\Vert_{2}\\
&  \leq3\,{\mathbb{E}}\big (|X^{(n)}_{i,j}|^{2}\mathbf{1}_{|X^{(n)}%
_{i,j}|>b_{n}}\big )+3\,{\mathbb{E}}\big (|X^{(n)}_{k,\ell}|^{2}%
\mathbf{1}_{|X^{(n)}_{k,\ell}|>b_{n}}\big )\,.
\end{align*}
Note also that, by assumption ${\mathbf{A}_{3}}$,
\begin{multline*}
\big |\mathrm{Cov}\big (X^{(n)}_{k,\ell}\mathbf{1}_{|X^{(n)}_{k,\ell}|\leq
b_{n}},X^{(n)}_{i,j}\mathbf{1}_{|X^{(n)}_{i,j}|\leq b_{n}}\big )-\mathrm{Cov}%
\big (X^{(n)}_{k,\ell},X^{(n)}_{i,j}\big )\big |\\
=\mathbf{1}_{i\in\lbrack k-K,k+K]}\mathbf{1}_{j\in\lbrack\ell-K,\ell
+K]}\big |\mathrm{Cov}\big (X^{(n)}_{k,\ell}\mathbf{1}_{|X^{(n)}_{k,\ell}|\leq
b_{n}},X^{(n)}_{i,j}\mathbf{1}_{|X^{(n)}_{i,j}|\leq b_{n}}\big )-\mathrm{Cov}%
\big (X^{(n)}_{k,\ell},X^{(n)}_{i,j}\big )\big |\,.
\end{multline*}
So, overall,
\begin{multline*}
\big |{\mathbb{E}}(G^{(n)}_{k,\ell}G^{(n)}_{i,j})-{\mathbb{E}}({\bar{G}}%
^{(n)}_{k,\ell}{\bar{G}}^{(n)}_{i,j})\big |\leq3\,{\mathbb{E}}\big (|X^{(n)}%
_{i,j}|^{2}\mathbf{1}_{|X^{(n)}_{i,j}|>b_{n}}\big )\mathbf{1}_{k\in\lbrack
i-K,i+K]}\mathbf{1}_{\ell\in\lbrack j-K,j+K]}\\
+3\,{\mathbb{E}}\big (|X^{(n)}_{k,\ell}|^{2}\mathbf{1}_{|X^{(n)}_{k,\ell
}|>b_{n}}\big )\mathbf{1}_{i\in\lbrack k-K,k+K]}\mathbf{1}_{j\in\lbrack
\ell-K,\ell+K]}\,.
\end{multline*}
Hence, starting from \eqref{approxtrungau} and taking into account the above
inequality, we derive that
\[
\big |{\mathbb{E}}\big (S_{{\bar{\mathbb{G}}}_{n}}(z)\big )-{\mathbb{E}%
}\big (S_{{\mathbb{G}}_{n}}(z)\big )\big |\leq\frac{12}{n^{2}v^{3}}%
(2K+1)^{2}\sum_{k=1}^{n}\sum_{\ell=1}^{k}{\mathbb{E}}\big (|X^{(n)}_{k,\ell
}|^{2}\mathbf{1}_{|X^{(n)}_{k,\ell}|>b_{n}}\big )\,.
\]
which converges to zero as $n$ tends to infinity, by assumption ${\mathbf{A}%
_{2}}$. This ends the proof of \eqref{aimthm1approxgau}.

\medskip

\noindent\textit{Step 2: Reduction to matrices with independent blocks.} Let
$p:=p_{n}$ such that $p_{n}\rightarrow\infty,$ $p_{n}/n\rightarrow0,$ and
$\tau_{n}p_{n}^{4}\rightarrow0.$ Clearly we can take $p_{n}>K$, $p_{n}+K\leq
n/3$, and set
\begin{equation}
q=q_{n}=\Big [\frac{n}{p+K}\Big ]-1\,. \label{defqk}%
\end{equation}
Let
\begin{equation}
I_{\ell}=\big [\ell(p+K)+1\,,\,\ell(p+K)+p\big ]\cap\mathbb{N}\ \text{
and}\ {\mathcal{E}}_{k,\ell}=\{(u,v)\in I_{k}\times I_{\ell-1}\}\,.
\label{defIE}%
\end{equation}
For any $k=1,\dots,q_{n}$ and any $\ell=1,\dots,k$, we define now real
matrices $B_{k,\ell}$ of size $p\times p$ whose entries consist of all the
${\bar{X}}_{u,v}^{(n)}$ for $(u,v)\in{\mathcal{E}}_{k,\ell}$. More precisely,
\begin{equation}
B_{k,\ell}:=\big \{b_{k,\ell}(i,j)\big \}_{i,j=1}^{p}\ \text{ where
$b_{k,\ell}(i,j)={\bar{X}}_{k(p+K)+i,(\ell-1)(p+K)+j}^{(n)}\,.$} \label{defBX}%
\end{equation}
Similarly, we define real matrices $B_{k,\ell}^{\ast}$ of size $p\times p$
whose entries consist of all the ${\bar{G}}_{u,v}^{(n)}$ for $(u,v)\in
{\mathcal{E}}_{k,\ell}$. Therefore,
\begin{equation}
B_{k,\ell}^{\ast}:=\big \{b_{k,\ell}^{\ast}(i,j)\big \}_{i,j=1}^{p}\ \text{
where $b_{k,\ell}^{\ast}(i,j)={\bar{G}}_{k(p+K)+i,(\ell-1)(p+K)+j}^{(n)}\,.$}
\label{defBY}%
\end{equation}
Using the blocks $B_{k,\ell}$ we construct now a $n\times n$ matrix,
${\widetilde{\mathbf{X}}}_{n},$ by inserting $0$'s. Actually we start from the
matrix ${\overline{\mathbf{X}}}_{n},$ keep the blocks $B_{k,\ell}$ and
$B_{k,\ell}^{T}$ and replace all the other variables by $0$'s. For the sake of
clarity we describe the south-western part of the matrix ${\widetilde
{\mathbf{X}}}_{n}$ below, the other part being constructed by symmetry.
\begin{equation}
{\widetilde{\mathbf{X}}}_{n}:=\left(
\begin{array}
[c]{ccccccccc}%
\mathbf{0}_{p,p} &  &  &  &  & ... &  &  & \\
\mathbf{0}_{K,p} & \mathbf{0}_{K,K} &  &  &  & ... &  &  & \\
B_{1,1} & \mathbf{0}_{p,K} & \mathbf{0}_{p,p} & \  &  & ... &  &  & \\
\mathbf{0}_{K,p} & \mathbf{0}_{K,K} & \mathbf{0}_{K,p} & \mathbf{0}_{K,K} &  &
... &  &  & \\
B_{2,1} & \mathbf{0}_{p,K} & B_{2,2} & \mathbf{0}_{p,K} & \mathbf{0}_{p,p} &
... &  &  & \\
... & ... & ... & ... & ... & ... & ... & ... & \mathbf{...}\\
B_{q-1,1} & \mathbf{0}_{p,K} & B_{q-1,2} & \mathbf{0}_{p,K} & B_{q-1,3} &
... &  &  & \\
\mathbf{0}_{K,p} & \mathbf{0}_{K,K} & \mathbf{0}_{K,p} & \mathbf{0}_{K,K} &
\mathbf{0}_{K,p} & ... & \mathbf{0}_{K,K} &  & \\
B_{q,1} & \mathbf{0}_{p,K} & B_{q,2} & \mathbf{0}_{p,K} & B_{q,3} & ... &
B_{q,q} & \mathbf{0}_{p,p} & \\
\mathbf{0}_{m,p} & \mathbf{0}_{m,K} & \mathbf{0}_{m,p} & \mathbf{0}_{m,K} &
\mathbf{0}_{m,p} & ... & \mathbf{0}_{m,K} & \mathbf{0}_{m,p} & \mathbf{0}%
_{m,m}%
\end{array}
\right)  \,, \label{defBX***}%
\end{equation}
where $m=m_{n}=n-q(p+K)-p$. ${\widetilde{\mathbf{G}}}_{n}$ is constructed as
${\widetilde{\mathbf{X}}}_{n}$ with the $B_{k,\ell}^{\ast}$ in place of the
$B_{k,\ell}$.

In what follows, we shall prove that, for any $z\in\mathbb{C}^{+}$,
\begin{equation}
\lim_{n\rightarrow\infty}\big |{\mathbb{E}}\big (S_{{\bar{\mathbb{X}}}_{n}%
}(z)\big )-{\mathbb{E}}\big (S_{{\widetilde{\mathbb{X}}}_{n}}%
(z)\big )\big |=0\ \text{ and }\ \lim_{n\rightarrow\infty}\big |{\mathbb{E}%
}\big (S_{{\bar{\mathbb{G}}}_{n}}(z)\big )-{\mathbb{E}}\big (S_{{\widetilde
{\mathbb{G}}}_{n}}(z)\big )\big |=0\,, \label{aimthm1inter}%
\end{equation}
with ${\widetilde{\mathbb{X}}}_{n}:=n^{-1/2}{\widetilde{\mathbf{X}}}_{n}$,
${\widetilde{\mathbb{G}}}_{n}:=n^{-1/2}{\widetilde{\mathbf{G}}}_{n}.$

To prove it we first introduce two other symmetric $n\times n$ matrices
${\widehat{\mathbf{X}}}_{n}=[{\widehat{X}}_{k,\ell}^{(n)}]_{k,\ell=1}^{n}$ and
${\widehat{\mathbf{G}}}_{n}=[{\widehat{G}}_{k,\ell}^{(n)}]_{k,\ell=1}^{n}$
constructed from ${\bar{\mathbf{X}}}_{n}$ and ${\bar{\mathbf{G}}}_{n}$
respectively, by replacing the entries by zeros in square blocks of size $p$
around the diagonal. More precisely, for any $1\leq i,j\leq n$,
\[
{\widehat{X}}_{i,j}^{(n)}=0\ \text{ if $(i,j)\in\cup_{\ell=0}^{q_{n}%
}{\mathcal{E}}_{\ell,\ell+1}$}\ \text{ and }\ {\widehat{X}}_{i,j}^{(n)}%
={\bar{X}}_{i,j}^{(n)}\ \text{ otherwise}%
\]
and
\[
{\widehat{G}}_{i,j}^{(n)}=0\ \text{ if $(i,j)\in\cup_{\ell=0}^{q_{n}%
}{\mathcal{E}}_{\ell,\ell+1}$}\ \text{ and }\ {\widehat{G}}_{i,j}^{(n)}%
={\bar{G}}_{i,j}^{(n)}\ \text{ otherwise}%
\]
where we recall that the sets ${\mathcal{E}}_{\ell,\ell+1}$ have been defined
in \eqref{defIE}. Denote now ${\widehat{\mathbb{X}}}_{n}=\frac{1}{\sqrt{n}%
}{\widehat{\mathbf{X}}}_{n}$ and ${\widehat{\mathbb{G}}}_{n}=\frac{1}{\sqrt
{n}}{\widehat{\mathbf{G}}}_{n}$.

By Lemma \ref{lmagotze}, we get, for any $z=u+\mathrm{i}v\in\mathbb{C}^{+}$,
that
\[
\big |{\mathbb{E}}\big (S_{{\bar{\mathbb{X}}}_{n}}(z)\big )-{\mathbb{E}%
}\big (S_{{\widehat{\mathbb{X}}}_{n}}(z)\big )\big |^{2}\leq{\mathbb{E}%
}\big (\big |S_{{\bar{\mathbb{X}}}_{n}}(z)-S_{{\widehat{\mathbb{X}}}_{n}%
}(z)\big |^{2}\big )\leq\frac{1}{n^{2}v^{4}}{\mathbb{E}}\big (\mathrm{Tr}%
\big (({\bar{\mathbf{X}}}_{n}-{\widehat{\mathbf{X}}}_{n})^{2}\big )\big )\,.
\]
Therefore,
\[
\big |{\mathbb{E}}\big (S_{{\bar{\mathbb{X}}}_{n}}(z)\big )-{\mathbb{E}%
}\big (S_{{\widehat{\mathbb{X}}}_{n}}(z)\big )\big |^{2}\leq\frac{1}%
{n^{2}v^{4}}\sum_{\ell=0}^{q_{n}}\sum_{(i,j)\in{\mathcal{E}}_{\ell,\ell+1}%
}{\mathbb{E}}(|{\bar{X}}_{i,j}^{(n)}|^{2})\,.
\]
But $\Vert{\bar{X}}_{i,j}^{(n)}\Vert_{\infty}\leq2\tau_{n}\sqrt{n}$. Hence,
\[
\big |{\mathbb{E}}\big (S_{{\bar{\mathbb{X}}}_{n}}(z)\big )-{\mathbb{E}%
}\big (S_{{\widehat{\mathbb{X}}}_{n}}(z)\big )\big |^{2}\leq\frac{4}%
{n^{2}v^{4}}(q_{n}+1)p_{n}^{2}\tau_{n}^{2}n\leq\frac{4}{v^{4}}\tau_{n}%
^{2}p_{n}\,.
\]
By our selection of $p_{n}$, we obviously have that $\tau_{n}^{2}%
p_{n}\rightarrow0$ as $n\rightarrow\infty$. It follows that, for any
$z\in\mathbb{C}^{+}$,
\begin{equation}
\lim_{n\rightarrow\infty}\big |{\mathbb{E}}\big (S_{{\bar{\mathbb{X}}}_{n}%
}(z)\big )-{\mathbb{E}}\big (S_{{\widehat{\mathbb{X}}}_{n}}%
(z)\big )\big |=0\,. \label{aimthm1interX*1}%
\end{equation}
With similar arguments, we get that, for any $z=u+\mathrm{i}v\in\mathbb{C}%
^{+}$,
\[
\big |{\mathbb{E}}\big (S_{{\bar{\mathbb{G}}}_{n}}(z)\big )-{\mathbb{E}%
}\big (S_{{\widehat{\mathbb{G}}}_{n}}(z)\big )\big |^{2}\leq\frac{1}%
{n^{2}v^{4}}\sum_{\ell=0}^{q_{n}}\sum_{(i,j)\in{\mathcal{E}}_{\ell,\ell}%
}{\mathbb{E}}(|{\bar{G}}_{i,j}^{(n)}|^{2})\,.
\]
But $\Vert{\bar{G}}_{i,j}^{(n)}\Vert_{2}=\Vert{\bar{X}}_{i,j}^{(n)}\Vert_{2}$.
So, as before, we derive that for any $z\in\mathbb{C}^{+}$,
\begin{equation}
\lim_{n\rightarrow\infty}\big |{\mathbb{E}}\big (S_{{\bar{\mathbb{G}}}_{n}%
}(z)\big )-{\mathbb{E}}\big (S_{{\widehat{\mathbb{G}}}_{n}}%
(z)\big )\big |=0\,. \label{aimthm1interG*1}%
\end{equation}
From \eqref{aimthm1interX*1} and \eqref{aimthm1interG*1}, the mean convergence
\eqref{aimthm1inter} follows if we prove that, for any $z\in\mathbb{C}^{+}$,
\begin{equation}
\lim_{n\rightarrow\infty}\big |{\mathbb{E}}\big (S_{{\widehat{\mathbb{X}}}%
_{n}}(z)\big )-{\mathbb{E}}\big (S_{{\widetilde{\mathbb{X}}}_{n}%
}(z)\big )\big |=0\ \text{ and }\ \lim_{n\rightarrow\infty}\big |{\mathbb{E}%
}\big (S_{{\widehat{\mathbb{G}}}_{n}}(z)\big )-{\mathbb{E}}%
\big (S_{{\widetilde{\mathbb{G}}}_{n}}(z)\big )\big |=0\,,
\label{aimthm1interbis}%
\end{equation}
For proving it, we shall use rank inequalities. Indeed, notice first that, for
any $z=u+\mathrm{i}v\in\mathbb{C}^{+}$,
\begin{align*}
\big |S_{{\widehat{\mathbb{X}}}_{n}}(z)  &  -S_{{\widetilde{\mathbb{X}}}_{n}%
}(z)\big |=\Big |\int\frac{1}{x-z}dF^{{\widehat{\mathbb{X}}}_{n}}(x)-\int
\frac{1}{x-z}dF^{{\widetilde{\mathbb{X}}}_{n}}(x)\Big |\\
&  \leq\Big |\int\frac{F^{{\widehat{\mathbb{X}}}_{n}}(x)-F^{{\widetilde
{\mathbb{X}}}_{n}}}{(x-z)^{2}}dx\Big |\leq\frac{\pi\,\big \|F^{{\widehat
{\mathbb{X}}}_{n}}-F^{{\widetilde{\mathbb{X}}}_{n}}\big \|_{\infty}}{v}\,.
\end{align*}
Hence, by Theorem A.43 in Bai and Silverstein \cite{BS},
\[
\big |S_{{\widehat{\mathbb{X}}}_{n}}(z)-S_{{\widetilde{\mathbb{X}}}_{n}%
}(z)\big |\leq\frac{\pi}{vn}\mathrm{rank}\big ({{\widehat{\mathbb{X}}}_{n}%
}-{{\widetilde{\mathbb{X}}}_{n}}\big )\,.
\]
But, by counting the numbers of rows and of columns with entries that can be
different from zero, we infer that
\[
\mathrm{rank}\big ({{\widehat{\mathbb{X}}}_{n}}-{{\widetilde{\mathbb{X}}}_{n}%
}\big )\leq2(q_{n}K+m_{n})\leq2(np^{-1}K+p+2K)\,.
\]
Therefore,
\[
\big |S_{{\widehat{\mathbb{X}}}_{n}}(z)-S_{{\widetilde{\mathbb{X}}}_{n}%
}(z)\big |\leq\frac{2\pi}{v}(Kp^{-1}+pn^{-1}+2Kn^{-1})\,.
\]
With similar arguments, we get
\[
\big |S_{{\widehat{\mathbb{G}}}_{n}}(z)-S_{{\widetilde{\mathbb{G}}}_{n}%
}(z)\big |\leq\frac{2\pi}{v}(Kp^{-1}+pn^{-1}+2Kn^{-1})\,.
\]
Since $p=p_{n}\rightarrow\infty$ and $p_{n}/n\rightarrow0$, as $n\rightarrow
\infty$, \eqref{aimthm1interbis} (and then \eqref{aimthm1inter}) follows from
the two above inequalities. Therefore, to prove that the mean convergence
\eqref{aimthm1} holds, it suffices to prove that, for any $z\in\mathbb{C}^{+}%
$,
\begin{equation}
\lim_{n\rightarrow\infty}\big |{\mathbb{E}}\big (S_{{\widetilde{\mathbb{X}}%
}_{n}}(z)\big )-{\mathbb{E}}\big (S_{{\widetilde{\mathbb{G}}}_{n}%
}(z)\big )\big |=0\,. \label{aimthm1*endstep2}%
\end{equation}
This is done in the next step.

\medskip

\noindent\textit{Step 3: Lindeberg method.} To prove \eqref{aimthm1*endstep2},
we shall use the Lindeberg method. Recall that the $\sigma$-fields $\sigma(
{\bar{G}}^{(n)}_{k,\ell} \, , \, 1 \leq\ell\leq k \leq n)$ and $\sigma(
X^{(n)}_{k,\ell} \, , \, 1 \leq\ell\leq k \leq n)$ are assumed to be
independent. Furthermore, by the hypothesis $\mathbf{A}_{3},$ all the blocks
$(B_{k,\ell})$ and $(B_{k,\ell}^{\ast})$ $1\leq\ell\leq k\leq q$ are independent.

The Lindeberg method consists of writing the difference of expectations as
telescoping sums and using Taylor expansions. This method can be used in the
context of random matrices since the function $f$, defined in \eqref{deffa},
admits partial derivatives of all orders (see the equality (17) in Chatterjee
\cite{Ct}). In the traditional Lindeberg method, the telescoping sums consist
of replacing one by one the random variables, involved in a partial sum, by a
Gaussian random variable. Here, we shall replace one by one the blocks
$B_{k,\ell}$ by the \textquotedblright Gaussian\textquotedblright\ ones
$B_{k,\ell}^{\ast}$ with the same covariance structure. So, starting from the
matrix ${\widetilde{\mathbf{X}}}_{n}={\widetilde{\mathbf{X}}}_{n}(0)$, the
first step is to replace its block $B_{q_{n},q_{n}}$ by $B_{q_{n},q_{n}}%
^{\ast}$, this gives a new matrix. Note that, at the same time, $B_{q_{n}%
,q_{n}}^{T}$ will also be replaced by $(B_{q_{n},q_{n}}^{\ast})^{T}.$ We
denote this matrix by ${\widetilde{\mathbf{X}}}_{n}(1)$ and re-denote the
block replaced by $B(1)$ and the new one by $B^{\ast}(1)$. At the second step,
we replace, in the new matrix ${\widetilde{\mathbf{X}}}_{n}(1)$, the block
$B(2):=B_{q_{n},q_{n}-1}$ by $B^{\ast}(2):=B_{q_{n},q_{n}-1}^{\ast}$, and call
the new matrix ${\widetilde{\mathbf{X}}}_{n}(2)$ and so on. Therefore, after
the $q_{n}$-th step, in the matrix ${\widetilde{\mathbf{X}}}_{n}$ we have
replaced the blocks $B(q_{n}-\ell+1)=B_{q_{n},\ell}\,,\,\ell=1,\dots,q_{n}$
(and their transposed) by the blocks $B^{\ast}(q_{n}-\ell+1)=B_{q_{n},\ell
}^{\ast}\,,\,\ell=1,\dots,q_{n}$ (and their transposed) respectively. This
matrix is denoted by ${\widetilde{\mathbf{X}}}_{n}(q_{n}).$ Next, the
$q_{n}+1$-th step will consist of replacing the block $B(q_{n}+1)=B_{q_{n}%
-1,q_{n}-1}$ by $B_{q_{n}-1,q_{n}-1}^{\ast}$ and obtain the matrix
${\widetilde{\mathbf{X}}}_{n}(q_{n}+1)$. So finally after $q_{n}(q_{n}+1)/2$
steps, we have replaced all the blocks $B_{k,\ell}\ $and $B_{k,\ell}^{T}$ of
the matrix ${\widetilde{\mathbb{X}}}_{n}$ to obtain at the end the matrix
${\widetilde{\mathbf{X}}}_{n}(q_{n}(q_{n}+1)/2)={\widetilde{\mathbb{G}}}_{n}$.

Therefore we have
\begin{equation}
{\mathbb{E}}\big (S_{{\widetilde{\mathbb{X}}}_{n}}(z)\big )-{\mathbb{E}%
}\big (S_{{\widetilde{\mathbb{G}}}_{n}}(z)\big )=\sum_{k=1}^{k_{n}%
}\Big( {\mathbb{E}}\big (S_{{\widetilde{\mathbb{X}}}_{n}(k-1)}%
(z)\big )-{\mathbb{E}}\big (S_{{\widetilde{\mathbb{X}}}_{n}(k)}%
(z)\big )\Big) \, . \label{telessum}%
\end{equation}
where $k_{n}=q_{n}(q_{n}+1)/2$.

Let $k$ in $\{1,\dots,k_{n}\}$. Observe that ${\widetilde{\mathbb{X}}}%
_{n}(k-1)$ and ${\widetilde{\mathbb{X}}}_{n}(k)$ differ only by the variables
in the block $B(k)$ replaced at the step $k$. Define then the vector
$\mathbf{X}$ of ${\mathbb{R}}^{p^{2}}$ consisting of all the entries of
$B(k)$, the vector $\mathbf{Y}$ of ${\mathbb{R}}^{p^{2}}$ consisting of all
the entries of $B^{\ast}(k)$ (in the same order we have defined the
coordinates of $\mathbf{X}$). Denote by $\mathbf{Z}$ the vector of
${\mathbb{R}}^{N-p^{2}}$ (where $N=n(n+1)/2$) consisting of all the entries on
and below the diagonal of ${\widetilde{\mathbf{X}}}_{n}(k-1)$ except the ones
that are in the block matrix $B(k)$. More precisely if $(u,v)$ are such that
$B(k)=B_{u,v}$, then
\[
\mathbf{X}=\big ((b_{u,v}(i,j))_{j=1,\dots,p}\,,\,i=1,\dots,p\big )\ \text{
and }\ \mathbf{Y}=\big ((b_{u,v}^{\ast}(i,j))_{j=1,\dots,p}\,,\,i=1,\dots
,p\big )
\]
where $b_{u,v}(i,j)$ and $b_{u,v}^{\ast}(i,j)$ are defined in \eqref{defBX}
and \eqref{defBY} respectively. In addition,
\[
\mathbf{Z}=\big (({\widetilde{\mathbf{X}}}_{n}(k-1))_{i,j}\,:\,1\leq j\leq
i\leq n\,,\,(i,j)\notin{\mathcal{E}}_{u,v}\big )\,,
\]
where ${\mathcal{E}}_{u,v}$ is defined in \eqref{defIE}. The notations above
allow to write
\[
{\mathbb{E}}\big (S_{{\widetilde{\mathbb{X}}}_{n}(k-1)}(z)\big )-{\mathbb{E}%
}\big (S_{{\widetilde{\mathbb{X}}}_{n}(k)}(z)\big )=\mathbb{E}f(\pi
(\mathbf{X},\mathbf{Z}))-\mathbb{E}f(\pi(\mathbf{Y},\mathbf{Z}))\,,
\]
where $f$ is the function from ${\mathbb{R}}^{N}$ to ${\mathbb{C}}$ defined by
\eqref{deffa} and $\pi:{\mathbb{R}}^{N}\rightarrow{\mathbb{R}}^{N}$ is a
certain permutation. Note that, by our hypothesis $\mathbf{A}_{3}$ and our
construction, the vectors $\mathbf{X}$, ${\mathbf{Y}}$ and $\mathbf{Z}$ are
independent. Moreover $\mathbf{X}$ and ${\mathbf{Y}}$ are centered at
expectation, have the same covariance structure and finite moments of order 3.
Applying then Lemma \ref{Taylor variables} from Section \ref{Sectappendix} and
taking into account \eqref{boundsd}, we derive that, for a constant $C$
depending only on $\operatorname{Im}(z)$,
\[
\big |{\mathbb{E}}\big (S_{{\widetilde{\mathbb{X}}}_{n}(k-1)}%
(z)\big )-{\mathbb{E}}\big (S_{{\widetilde{\mathbb{X}}}_{n}(k)}%
\big )\big |\leq\frac{Cp^{4}}{n^{5/2}}\sum_{(i,j)\in{\mathcal{E}}_{u,v}%
}\left(  {\mathbb{E}}(|\bar{X}_{i,j}^{(n)}|^{3})+{\mathbb{E}}(|\bar{G}%
_{i,j}^{(n)}|^{3})\right)  \,.
\]
So, overall,
\begin{equation}
\sum_{k=1}^{k_{n}}|{\mathbb{E}}\big (S_{{\widetilde{\mathbb{X}}}_{n}%
(k-1)}(z)\big )-{\mathbb{E}}\big (S_{{\widetilde{\mathbb{X}}}_{n}%
(k)}\big )|\leq\frac{Cp^{4}}{n^{5/2}}\sum_{1\leq\ell\leq k\leq q}%
\sum_{(i,j)\in{\mathcal{E}}_{k,\ell}}\left(  {\mathbb{E}}(|\bar{X}_{i,j}%
^{(n)}|^{3})+{\mathbb{E}}(|\bar{G}_{i,j}^{(n)}|^{3})\right)  \,.
\label{estimate 1}%
\end{equation}
By taking into account that
\[
{\mathbb{E}}(|\bar{X}_{i,j}^{(n)}|^{3})\leq2\tau_{n}\sqrt{n}{\mathbb{E}%
}(|X_{i,j}^{(n)}|^{2})
\]
and also
\[
{\mathbb{E}}(|\bar{G}_{i,j}^{(n)}|^{3})\leq2\,\big ({\mathbb{E}}(|\bar
{G}_{i,j}^{(n)}|^{2})\big )^{3/2}=2\,\big ({\mathbb{E}}(|\bar{X}_{i,j}%
^{(n)}|^{2})\big )^{3/2}\leq4\tau_{n}\sqrt{n}{\mathbb{E}}(|X_{i,j}^{(n)}%
|^{2})  \,,
\]
it follows from (\ref{telessum}) and (\ref{estimate 1}) that, for a constant
$C^{\prime}$ depending only on $\operatorname{Im}(z)$,
\[
\big |{\mathbb{E}}\big (S_{{\widetilde{\mathbb{X}}}_{n}}(z)\big )-{\mathbb{E}%
}\big (S_{{\widetilde{\mathbb{G}}}_{n}}(z)\big ) \big |\leq\frac{C^{\prime
}p^{4}}{n^{2}}\tau_{n}\sum_{1\leq j\leq i\leq n}{\mathbb{E}}(|X_{i,j}%
^{(n)}|^{2})
\]
which converges to $0$ by $\mathbf{A}_{1}$ and the selection of $p_{n}$. This
ends the proof of \eqref{aimthm1*endstep2} and then of the theorem.
\hfill$\square$

\medskip

\noindent\textbf{Acknowledgements.} The authors would like to thank the
referee for carefully reading the manuscript and J. Najim for helpful discussions.


\begin{thebibliography}{99}                                                                                               %
\bibitem {AZ}Anderson, G. and Zeitouni, O. (2008). A law of large numbers for
finite-range dependent random matrices \textit{Comm. Pure Appl. Math.}
\textbf{61} 1118-1154.

\bibitem {BS}Bai, Z. and Silverstein, J.W. (2010). \textit{Spectral analysis
of large dimensional random matrices}. Springer, New York, second edition.

\bibitem {BZ}Bai, Z. and Zhou, W. (2008). Large sample covariance matrices
without independence structures in columns. \emph{Statist. Sinica} \textbf{18} 425-442.

\bibitem {BM}Banna, M. and Merlev\`{e}de, F. (2013). Limiting spectral
distribution of large sample covariance matrices associated with a class of
stationary processes. To appear in \textit{J. Theoret. Probab.} (DOI: 10.1007/s10959-013-0508-x)

\bibitem {BK}Boutet de Monvel, A. and Khorunzhy, A. (1999). On the Norm and
Eigenvalue Distribution of Large Random Matrices. \textit{ Ann. Probab.}
\textbf{27} 913-944.

\bibitem {BKV}Boutet de Monvel, A. Khorunzhy, A. and Vasilchuk, V. (1996).
Limiting eigenvalue distribution of random matrices with correlated entries.
\textit{Markov Process. Related Fields} \textbf{2} 607-636.

\bibitem {chak}Chakrabarty A., Hazra R.S. and Sarkar D. (2014). From random
matrices to long range dependence. \textit{arXiv:math/1401.0780}.

\bibitem {Ct}Chatterjee, S. (2006). A generalization of the Lindeberg
principle.\textit{ Ann. Probab.} \textbf{34} 2061-2076.





\bibitem {Girko90}Girko, V. L. (1990). Theory of Random Determinants.
Translated from the Russian. \textit{Mathematics and Its Applications (Soviet
Series)} \textbf{45}. Kluwer Academic Publishers Group, Dordrecht.

\bibitem {GNT}G\"{o}tze, F., Naumov, A. and A. Tikhomirov (2012). Semicircle
law for a class of random matrixes with dependent entries. arXiv:math/0702386v1.

\bibitem {GL}Guntuboyina, A. and Leeb, H. (2009). Concentration of the
spectral measure of large Wishart matrices with dependent entries.
\textit{Electron. Commun. Probab.} \textbf{14} 334-342.

\bibitem {HLN}Hachem, W., Loubaton, P. and J. Najim (2005). The empirical
eigenvalue distribution of a Gram matrix: from independence to stationarity,
\textit{Markov Process. Related Fields} \textbf{11} 629--648.



\bibitem {KP}Khorunzhy, A. and Pastur, L. (1994). On the eigenvalue
distribution of the deformed Wigner ensemble of random matrices. In: V. A.
Marchenko (ed.), \textit{Spectral Operator Theory and Related Topics}, Adv.
Soviet Math. \textbf{19}, Amer. Math. Soc., Providence, RI, 97-127.

\bibitem {MP}Merlev\`ede, F. and Peligrad, M. (2014). On the empirical
spectral distribution for matrices with long memory and independent rows.
\textit{arXiv: 1406.1216}

\bibitem {PS}Pastur, L. and Shcherbina, M. (2011). Eigenvalue distribution of
large random matrices. \textit{Mathematical Surveys and Monographs},
\textbf{171}. American Mathematical Society, Providence, RI.



\bibitem {ROBS}Rashidi Far, R., Oraby T., Bryc, W. and Speicher, R. (2008). On
slow-fading MIMO systems with nonseparable correlation. \textit{IEEE Trans.
Inform. Theory} \textbf{54} 544-553.

\bibitem {Tal}Talagrand M. (2010). \textit{Mean Field Models for Spin Glasses.
Vol 1. Basic Examples.} Springer.



\bibitem {Yao}Yao, J. (2012). A note on a Mar\u{c}enko-Pastur type theorem for
time series. \textit{Statist. Probab. Lett.} \textbf{82} 22-28.
\end{thebibliography}
\end{document}